\documentstyle[11pt,amssymb,amsfonts]{article}

\begin{document}
\newcommand{\p}{\parallel }
\makeatletter \makeatother
\newtheorem{th}{Theorem}[section]
\newtheorem{lem}{Lemma}[section]
\newtheorem{de}{Definition}[section]
\newtheorem{rem}{Remark}[section]
\newtheorem{cor}{Corollary}[section]
\renewcommand{\theequation}{\thesection.\arabic {equation}}

\title{{\bf Multiply twisted products}
}
\author{Yong Wang \\}

\date{}
\maketitle

\begin{abstract} In this paper, we compute the index form of the multiply twisted products. We study the Killing vector fields
 on the multiply twisted product manifolds and determine the  Killing vector
 fields in some cases. We compute the curvature of the multiply twisted
 products with a semi-symmetric metric connection and show that the
 mixed Ricci-flat multiply twisted products with a semi-symmetric metric
 connection can be expressed as multiply warped products. We also
 study the Einstein multiply warped products with a semi-symmetric metric
 connection and the multiply warped products with a semi-symmetric metric
 connection with constant scalar curvature, we apply our results to
 generalized Robertson-Walker spacetimes with a semi-symmetric metric connection
 and generalized Kasner spacetimes with a semi-symmetric metric
 connection and find some new examples of Einstein affine
 manifolds and affine manifolds with constant scalar curvature.
  We also consider the multiply twisted product Finsler
 manifolds and we get some interesting properties of these spaces.\\

\end{abstract}

\section{Introduction}
    \quad  The (singly) warped product $B\times_bF$  of two pseudo-Riemannian manifolds $(B,g_B)$ and
    $(F,g_F)$ with a smooth function $b:B\rightarrow (0,\infty)$ is the product
    manifold $B\times F$ with the metric tensor $g=g_B\oplus
    b^2g_F.$ Here, $(B,g_B)$ is called the base manifold and
    $(F,g_F)$ is called as the fiber manifold and $b$ is called as
    the warping function. Generalized Robertson-Walker space-times
    and standard static space-times are two well-known warped
    product spaces.The concept of warped products was first introduced by
    Bishop and ONeil (see [BO]) to construct examples of Riemannian
    manifolds with negative curvature. In Riemannian geometry,
    warped product manifolds and their generic forms have been used
    to construct new examples with interesting curvature properties
    since then. In [DD], F. Dobarro and E. Dozo had studied from the viewpoint of partial differential equations and variational methods,
    the problem of showing when a Riemannian metric of constant scalar curvature can be produced on a product manifolds by a warped product construction.
    In [EJK], Ehrlich, Jung and Kim got explicit solutions to warping function to have a constant scalar curvature for generalized Robertson-Walker space-times.
    In [ARS], explicit solutions were also obtained for the warping
    function to make the space-time as Einstein when the fiber is
    also Einstein.\\
      \indent One can generalize singly warped products to multiply warped
      products. Briefly, a multiply warped product $(M,g)$ is a
      product manifold of form
      $M=B\times_{b_1}F_1\times_{b_2}F_2\cdots \times_{b_m}F_m$ with
      the metric $g=g_B\oplus b_1^2g_{F_1}\oplus
      b_2^2g_{F_2}\cdots \otimes b_m^2g_{F_m}$, where for each $i\in
      \{ 1,\cdots,m\},~b_i:~B\rightarrow (0,\infty)$ is smooth and
      $(F_i,g_{F_i})$ is a pseudo-Riemannian manifold. In
      particular, when $B=(c,d)$ with the negative definite metric
      $g_B=-dt^2$ and $(F_i,g_{F_i})$ is a Riemannian manifold, we
      call $M$ as the multiply generalized Robertson-Walker
      space-time. Geodesic equations and geodesic connectedness of
      multiply generalized Robertson-Walker
      space-times were studied by Flores and S\'{a}nchez in [FS] and
      they also noted that the class of multiply generalized warped
      space-times contains many well known relativistic space-times.
      In [U1], necessary and sufficient conditions were obtained
      about geodesic completeness of multiply warped space-times. In
      [DU1], Dobarro and \"{U}nal studied Ricci-flat and
      Einstein-Lorentzian multiply warped products and considered
      the case of having constant scalar curvature for multiply warped
      products and applied their results to generalized Kasner
      space-times.\\
      \indent Singly warped products have two natural generalizations. A doubly warped product
      $(M,g)$ is a
      product manifold of form
      $M=_fB\times_{b}F$, with smooth functions $b:B\rightarrow (0,\infty)$, $f:F\rightarrow (0,\infty)$
      and the metric tensor $g=f^2g_B\oplus b^2g_F.$ In [U2], \"{U}nal studied geodesic completeness of Riemannian
      doubly warped products and Lorentzian doubly warped products. A twisted product
      $(M,g)$ is a
      product manifold of form
      $M=B\times_{b}F$, with a smooth function $b:B\times F\rightarrow (0,\infty)$,
      and the metric tensor $g=g_B\oplus b^2g_F.$ In [FGKU], they
      showed that mixed Ricci-flat twisted products could be expressed as warped
      products. As a consequence, any Einstein twisted products are
      warped products. So it is natural to consider multiply twisted
      products as generalizations of multiply warped
      products and twisted products. A multiply twisted product $(M,g)$ is a
      product manifold of form
      $M=B\times_{b_1}F_1\times_{b_2}F_2\cdots \times_{b_m}F_m$ with
      the metric $g=g_B\oplus b_1^2g_{F_1}\oplus
      b_2^2g_{F_2}\cdots \oplus b_m^2g_{F_m}$, where for each $i\in
      \{ 1,\cdots,m\},~b_i:~B\times F_i\rightarrow (0,\infty)$ is
      smooth.\\
      \indent In [EK], Ehrlich and Kim constructed the index form
      along timelike geodesics on a Lorentzian warped product and
      applied this index form to generalized Robertson-Walker
      spacetimes. In [CK], they computed the index form of multiply generalized Robertson-Walker
      spacetimes. In the first part of this paper, we compute the index form
      of multiply twisted products. In [Sa], the curvature and
      Killing vector fields of generalized Robertson-Walker
      spacetimes were studied. The non-trivial Killing vector fields
      of these space-times were characterized. In [DU2], they
      provided a global characterization of the Killing vector
      fields of a standard static spacetime by a system of partial
      differential equations. By studying this system, they
      determined all the Killing vector fields when Riemannian part
      was compact. In the second part of this paper, we study the Killing vector fields
        of multiply twisted products and in some cases, we can determine all Killing vector fields.\\
        \indent The definition of a semi-symmetric metric connection was given by H. Hayden in [Ha]. In 1970, K. Yano [Ya]
        considered a semi-symmetric metric connection and studied some of its properties. He proved that a Riemannian manifold admitting
         the semi-symmetric metric connection has vanishing curvature tensor if and only if it is conformally flat. Motivated by the Yano' result,
         in [SO], Sular and \"{O}zgur studied warped product manifolds with a
         semi-symmetric metric connection, they computed curvature of semi-symmetric metric connection
          and considered Einstein warped product manifolds with a semi-symmetric metric connection. In the main part of this paper,
           we consider multiply twisted products with a semi-symmetric metric connection and compute the curvature of a semi-symmetric metric connection.
            We
         show that mixed Ricci-flat multiply twisted products with a semi-symmetric metric
 connection can be expressed as multiply warped products which
 generalizes the result in [FGKU]. We also
 study the Einstein multiply warped products with a semi-symmetric metric
 connection and multiply warped products with a semi-symmetric metric
 connection with constant scalar curvature, we apply our results to
 generalized Robertson-Walker spacetimes with a semi-symmetric metric connection
 and generalized Kasner spacetimes with a semi-symmetric metric
 connection and we find some new examples of Einstein affine
 manifolds and affine manifolds with constant scalar curvature. We
 also classify generalized Einstein Robertson-Walker spacetimes with a semi-symmetric metric
 connection and generalized Einstein Kasner spacetimes with a semi-symmetric metric
 connection.\\
      \indent On the other hand, Finsler geometry is a subject
      studying manifolds whose tangent spaces carry a norm varying
      smoothly with the base point. Indeed, Finsler geometry is just
      Riemannian geometry without the quadratic restriction. Thus it
      is natural to extend the construction of warped product
      manifolds for Finsler geometry. In the first step, Asanov gave
      the generalization of the Schwarzschild metric in the
      Finslerian setting and obtained some models of relativety theory
      descried through the warped product of Finsler metrics
      [As1,2]. In [KPV], Kozma-Peter-Varga defined their warped
      product for Finsler metrics and concluded that completeness of
      warped product can be related to completeness of its
      components. In [HR], Using the Cantan connection for the study
      of warped product Finsler spaces, they found the necessary and
      sufficient conditions for such manifolds to be Riemannian,
      Landsberg, Berwald, and locally Minkowski, separately. In [PT]
      and [PTN], they considered the doubly warped product Finsler
      manifolds and found the necessary and
      sufficient conditions for such manifolds to be Riemannian,
      Landsberg, Berwald, Douglas, locally dually flat. They also
      defined the doubly warped Sasaki-Matsumoto metric for warped
      product manifolds and found a condition under which the
      horizontal and vertical tangent bundle were totally geodesic.
      In this paper, we consider multiply twisted product Finsler
      manifolds. Let $(B,\overline{F_B}),~(M_i,\overline{F_i}),~1\leq i \leq m$ be Finsler manifolds and
      $f_i: B\times M_i\rightarrow (0,+\infty)$ be a
      smooth function. Let $\pi _i: TM_i\rightarrow M_i$ be the projection map. The product manifold
      $B\times M_1\times \cdots \times M_m$
      endowed with the metric $\overline{F}:TB^0\times TM_1^0\cdots\times
      TM_m^0\rightarrow \mathbb{R}$ is considered,
      $$\overline{F}(v_0,v_1,\cdots,v_m)$$
      $$=\sqrt{\overline{F_B}^2(v_0)+f_1^2(\pi_0(v_0),\pi_1(v_1))\overline{F_1}^2(v_1)+\cdots+
     f_m^2(\pi_0(v_0),\pi_m(v_m))\overline{F_m}^2(v_m)},\eqno(1.1)$$
     where $TB^0=TB-\{0\},~TM^0_i=TM_i-\{0\}.$ We find the necessary and
      sufficient conditions for multiply twisted product Finsler
      manifolds to be Riemannian, Landsberg, Berwald, locally dually flat, locally Minkowski.We also
      define the generalized twisted Sasaki-Matsumoto metric and get a condition under which the
      horizontal and vertical tangent bundle are totally
      geodesic.\\
      \indent This paper is arranged as follows: In Section 2,  we compute curvature and the index form
      of multiply twisted products. In Section 3,  we study the Killing vector fields
        of multiply twisted products and in some cases, we can determine all Killing vector
        fields. In Section 4, we study multiply twisted products with a semi-symmetric metric
        connection. In Section 5, we study multiply twisted product Finsler
      manifolds.\\

\section{ The index form of the multiply twisted products}

\noindent {\bf Definition 2.1}  A {\it multiply twisted product}
$(M,g)$ is a product manifold of form
      $M=B\times_{b_1}F_1\times_{b_2}F_2\cdots \times_{b_m}F_m$ with
      the metric $g=g_B\oplus b_1^2g_{F_1}\oplus
      b_2^2g_{F_2}\cdots \oplus b_m^2g_{F_m}$, where for each $i\in
      \{ 1,\cdots,m\},~b_i:~B\times F_i\rightarrow (0,\infty)$ is
      smooth.\\
\indent Here, $(B,g_B)$ is called the base manifold and
    $(F_i,g_{F_i})$ is called as the fiber manifold and $b_i$ is called as
    the twisted function. Obviously, twisted products and multiply warped
    products are the special cases of multiply twisted products.\\

    \noindent {\bf Proposition 2.2} {\it Let $M=B\times_{b_1}F_1\times_{b_2}F_2\cdots
    \times_{b_m}F_m$ be a multiply twisted product and let $X,Y\in \Gamma(TB)$
      and $U\in \Gamma(TF_i)$, $W\in \Gamma(TF_j)$. Then}\\
      \noindent$(1) ~~\nabla_XY=\nabla^B_XY.$\\
      \noindent$(2)~~\nabla_XU=\nabla_UX=\frac{X(b_i)}{b_i}U.$\\
            \noindent$(3)~~\nabla_UW=0~~if~ i\neq j.$\\
      \noindent$(4)~~\nabla_UW=U({\rm ln}b_i)W+W({\rm ln}b_i)U-\frac{g_{F_i}(U,W)}{b_i}{\rm
      grad}_{F_i}b_i-b_ig_{F_i}(U,W){\rm
      grad}_{B}b_i+\nabla^{F_i}_UW~~if~ i= j.$\\

      \noindent {\bf Proposition 2.3} {\it Let $M=B\times_{b_1}F_1\times_{b_2}F_2\cdots
    \times_{b_m}F_m$ be a multiply twisted product, then $B$ is a
    totally geodesic submanifold and $F_i$ is a totally umbilical
    submanifold.}\\

    Define the curvature, Ricci curvature and scalar curvature as
    follows:
    $$R(X,Y)Z=\nabla_X\nabla_Y-\nabla_Y\nabla_X-\nabla_{[X,Y]},$$
    $$ ~{\rm
    Ric}(X,Y)=\sum_k\varepsilon_k<R(X,E_k)Y,E_k>, ~~S=\sum_k\varepsilon_k{\rm
    Ric}(E_k,E_k),$$ where $E_k$ is a orthonormal base of $M$ with
    $<E_k,E_k>=\varepsilon_k.$ The Hessian of $f$ is defined by
    $H^f(X,Y)=XYf-(\nabla_XY)f.$\\

     \noindent {\bf Proposition 2.4} {\it Let $M=B\times_{b_1}F_1\times_{b_2}F_2\cdots
    \times_{b_m}F_m$ be a multiply twisted product and let $X,Y,Z\in \Gamma(TB)$
      and $V\in \Gamma(TF_i)$, $W\in \Gamma(TF_j)$, $U\in \Gamma(TF_k)$. Then}\\
\noindent $(1)R(X,Y)Z=R^B(X,Y)Z.$\\
\noindent $(2)R(V,X)Y=-\frac{H^{b_i}_B(X,Y)}{b_i}V.$\\
\noindent $(3)R(X,V)W=R(V,W)X=R(V,X)W=0 ~if ~i\neq j.$\\
\noindent $(4)R(X,Y)V=0.$\\
\noindent $(5)R(V,W)X=VX({\rm ln}b_i)W-WX({\rm ln}b_i)V ~ if ~i=j.$\\
\noindent $(6)R(V,W)U=0~ if ~i=j\neq k~ or~ i\neq j \neq k.$\\
\noindent $(7)R(U,V)W=-g(V,W)\frac{g_B({\rm grad}_Bb_i,{\rm
grad}_Bb_k)}{b_ib_k}U, ~if ~i=j\neq k.$\\
\noindent $(8)R(X,V)W=-\frac{g(V,W)}{b_i}\nabla_X^B({\rm
grad}_Bb_i)+[WX({\rm ln}b_i)]V-g_{F_i}(W,V){\rm grad}_{F_i}(X{\rm
ln}b_i) ~if
~i=j.$\\
\noindent $(9)R(V,W)U=g(V,U){\rm grad}_B(W({\rm ln}b_i))-g(W,U){\rm
grad}_B(V({\rm ln}b_i))+R^{F_i}(V,W)U-\frac{|{\rm
grad}_Bb_i|^2_B}{b_i^2}(g(W,U)V-g(V,U)W) ~if ~i=j=k.~$\\

\noindent {\bf Proposition 2.5} {\it Let
$M=B\times_{b_1}F_1\times_{b_2}F_2\cdots
    \times_{b_m}F_m$ be a multiply twisted product and let $X,Y,Z\in \Gamma(TB)$
      and $V\in \Gamma(TF_i)$, $W\in \Gamma(TF_j)$. Then}\\
\noindent $(1) {\rm Ric} (X,Y)={\rm
Ric}^B(X,Y)+\sum_{i=1}^m\frac{l_i}{b_i}H_B^{b_i}(X,Y).$\\
\noindent $(2) {\rm Ric} (X,V)={\rm Ric}
(V,X)=(l_i-1)[VX({\rm ln}b_i)].$\\
\noindent $(3) {\rm Ric} (V,W)=0~if ~i\neq j.$\\
\noindent $(4) {\rm Ric} (V,W)={\rm Ric}^{F_i}
(V,W)+\left[\frac{\triangle_Bb_i}{b_i}+(l_i-1)\frac{|{\rm
grad}_Bb_i|^2_B}{b_i^2}+\sum_{k\neq i}l_k\frac{g_B({\rm
grad}_Bb_i,{\rm grad}_Bb_k)}{b_ib_k}\right]$\\
$~~~~~~\cdot g(V,W)~if ~i= j.$\\

\indent By Proposition 2.5, similar to the theorem 1 in [FGKU], we
get:\\

 \noindent {\bf Corollary 2.6} {\it Let
$M=B\times_{b_1}F_1\times_{b_2}F_2\cdots
    \times_{b_m}F_m$ be a multiply twisted product and ${\rm
    dim}F_i>1$, then $M$ is mixed Ricci-flat if and only if $M$ can
    be expressed as a multiply warped product. In particular, if $M$
    is Einstein, then $M$ can
    be expressed as a multiply warped product.}\\

\indent Similar to the theorem 6 in [BGV], we
get:\\

 \noindent {\bf Corollary 2.7} {\it Let
$M=B\times_{b_1}F_1\times_{b_2}F_2\cdots
    \times_{b_m}F_m$ be a multiply twisted product and ${\rm
    dim}B>1$, ${\rm
    dim}F_i>1$, if $M$ is locally conformal-flat, then $M$ can
    be expressed as a multiply warped product.}\\

\noindent {\bf Proposition 2.8} {\it Let
$M=B\times_{b_1}F_1\times_{b_2}F_2\cdots
    \times_{b_m}F_m$ be a multiply twisted product, then the scalar
    curvature $S$ has the following expression:}\\
    $$
    S=S^B+2\sum_{i=1}^m\frac{l_i}{b_i}\triangle_Bb_i+\sum_{i=1}^m\frac{S^{F_i}}{b_i^2}+\sum_{i=1}^ml_i(l_i-1)
\frac{|{\rm grad}_Bb_i|^2_B}{b_i^2}+\sum_{i=1}^m\sum_{k\neq
i}l_il_k\frac{g_B({\rm grad}_Bb_i,{\rm
grad}_Bb_k)}{b_ib_k}.\eqno(2.1)$$

\indent Similar to the proposition 3.1 in [U2], we have\\

\noindent {\bf Proposition 2.9} {\it Let
$M=B\times_{b_1}F_1\times_{b_2}F_2\cdots
    \times_{b_m}F_m$ be a multiply twisted product. If $(B,g_B)$ and
    $(F_i,g_{F_i})$ are complete Riemannian manifolds and ${\rm inf}
    b_i>0$, then $(M,g)$ is a complete Riemannian manifold.
    Conversely, if  $(M,g)$ is a complete Riemannian manifold and
    $0<{\rm inf}
    b_i<{\rm sup}b_i<+\infty$, then $(B,g_B)$ and
    $(F_i,g_{F_i})$ are complete Riemannian manifolds.}\\

\indent By Proposition 2.2, we have:\\

\noindent {\bf Proposition 2.10} {\it Let
$M=B\times_{b_1}F_1\times_{b_2}F_2\cdots
    \times_{b_m}F_m$ be a multiply twisted product. Also let
    $\gamma=(\alpha,\beta_1,\cdots,\beta_m)$ be a curve in $M$
    defined on some interval $I\subseteq {\bf R}.$ Then $\gamma$ is
    a geodesic in $M$ if and only if for any $t\in I$,}\\
    \indent
    $1.\alpha^{''}=\sum_{i=1}^mb_i(\alpha,\beta_i)g_{F_i}(\beta_i',\beta_i'){\rm
    grad}_Bb_i.$\\
\indent
    $2.\beta_i''=\frac{-2\alpha'(b_i)}{b_i(\alpha,\beta_i)}\beta_i'-2\beta_i'(lnb_i)\beta_i'+\frac{g_{F_i}(\beta_i',\beta_i')}{b_i(\alpha,\beta_i)}
{\rm grad}_{F_i}b_i,~for~ any ~i\in\{1,2\cdots,m\}.$\\

\indent By the proposition 2.2, the formula for the covariant
derivative of a smooth vector field $V=(V_B,V_{F_1},\cdots,V_{F_m})$
along the smooth curve $\gamma=(\alpha,\beta_1,\cdots,\beta_m)$ may
be obtained:
$$V'(t)=\left(V_B'(t)-\sum_{i=1}^mb_i<\beta_i'(t),V_i>_{F_i}{\rm
grad}_Bb_i\right.,$$
$$\left.\sum_{i=1}^m[\frac{\alpha'(b_i)}{b_i(\alpha,\beta_i)}V_i+\frac{V_B(b_i)}{b_i}\beta'_i(t)+\beta_i'({\rm ln}b_i)V_i+V_i({\rm ln}b_i)\beta'_i
-\frac{g_{F_i}(\beta'_i,V_i)}{b_i}{\rm
grad}_{F_i}b_i+\nabla^{F_i}_{\beta'_i}V_i]\right).\eqno(2.2)$$

\indent By (2.2), we have:\\

\noindent {\bf Proposition 2.11} {\it Let
$M=B\times_{b_1}F_1\times_{b_2}F_2\cdots
    \times_{b_m}F_m$ be a multiply twisted product.\\
\noindent (a) Let $V_B$ be a smooth vector field along the smooth
curve $\gamma_B:I\rightarrow (B, g_B)$ and for fixed $q_i\in F_i$,
let $\gamma(t)=(\gamma_B(t),q_1,\cdots,q_m)$ and
$V(t)=(V_B(t),0,\cdots,0)$ along $\gamma$. Then $V$ is parallel
along $\gamma:I\rightarrow (M,g)$ if and only if $V_B$ is parallel
along $\gamma_B:I\rightarrow (B,g_B).$\\
\noindent (b) Let $V_{F_i}$ be a smooth vector field along the
smooth curve $\gamma_{F_i}:I\rightarrow (F_i, g_{F_i})$ and for
fixed $b\in B$ and $q_j\in F_j,~j\neq i$, let
$\gamma(t)=(b,q_1,\cdots,q_{i-1},\gamma_{F_i}(t),\cdots,q_m)$ and
$V(t)=(0,\cdots,0,V_{F_i}(t),\cdots,0)$ along $\gamma$. Then $V$ is
parallel along $\gamma:I\rightarrow (M,g)$ if and only if
$<\gamma'_{F_i}(t),V_{F_i}(t)>_{F_i}{\rm grad}_Bb_i=0$ and}
$$\gamma'_{F_i}(t)(lnb_i)V_i+V_i(lnb_i)\gamma'_{F_i}(t)
-\frac{g_{F_i}(\gamma'_{F_i}(t),V_i)}{b_i}{\rm
grad}_{F_i}b_i+\nabla^{F_i}_{\gamma'_{F_i}(t)}V_i=0.$$\\

\indent Next we compute the index form. We use the variational
approach rather than direct computation of the curvature formula
$g(R(V,\gamma')\gamma',V)$ of the index form.\\
\indent Now let $\gamma: [a,b]\rightarrow (M,g)$ be a unit timelike
curve. Further, let $\alpha:I\times
(-\varepsilon,\varepsilon)\rightarrow(M,g)$ be a variation of
$\gamma(t)$. Let
$$\alpha_s(t):=\alpha(t,s)=(\alpha_B(t,s),\alpha_1(t,s),\cdots,\alpha_m(t,s))\eqno(2.3)$$
and define corresponding variation vector fields
$$W=\alpha_*\frac{\partial}{\partial s}=\frac{\partial \alpha}{\partial
s};~W_B=(\alpha_B)_*\frac{\partial}{\partial s}=\frac{\partial
\alpha_B}{\partial s};W_i=(\alpha_i)_*\frac{\partial}{\partial
s}=\frac{\partial \alpha_i}{\partial s}\eqno(2.4)$$ \noindent
$$V(t)=W(t,0),~~ V_B(t)=W_B(t,0),~~V_i(t)=W_i(t,0)\eqno(2.5)$$

\noindent Also,
$$\alpha'_s=(\alpha_{B_*}\frac{\partial}{\partial t},{\alpha_1}_*\frac{\partial}{\partial
t},\cdots,{\alpha_m}_*\frac{\partial}{\partial t})=(\frac{\partial
\alpha_B}{\partial t},\frac{\partial \alpha_1}{\partial
t},\cdots,\frac{\partial \alpha_m}{\partial t}).\eqno(2.6)$$ Since
the curve $\alpha_s(t)$ is timelike, if
$$h(t,s)=-g_B(\frac{\partial
\alpha_B}{\partial t},\frac{\partial \alpha_B}{\partial t})-\sum
_{i=1}^mb_i^2(\alpha_B,\alpha_i)(t,s)g_{F_i}(\frac{\partial
\alpha_i}{\partial t},\frac{\partial \alpha_i}{\partial
t})=-g(\frac{\partial \alpha}{\partial t},\frac{\partial
\alpha}{\partial t})\eqno(2.7)$$ \noindent and
$$A(t,s)=g_B(\nabla^B_{\frac{\partial }{\partial s}}\frac{\partial
\alpha_B}{\partial t},\frac{\partial \alpha_B}{\partial
t})+\frac{1}{2}\sum _{i=1}^m\frac{\partial }{\partial
s}(b_i^2(\alpha_B,\alpha_i))(t,s)g_{F_i}(\frac{\partial
\alpha_i}{\partial t},\frac{\partial \alpha_i}{\partial t})$$ $$ +
\sum
_{i=1}^mb_i^2(\alpha_B,\alpha_i)(t,s)g_{F_i}(\nabla_{\frac{\partial
}{\partial s}}^{F_i}\frac{\partial \alpha_i}{\partial
t},\frac{\partial \alpha_i}{\partial t})\eqno(2.8)$$ \noindent then
$$h(t,0)=1,~~~\frac{\partial h}{\partial s}=-2A(t,s)\eqno(2.9)$$
\noindent and
$$L(\alpha_s)=\int_{t=a}^b\sqrt{-g(\frac{\partial \alpha}{\partial t},\frac{\partial
\alpha}{\partial
t})}=\int_{t=a}^b(h(t,s)^{\frac{1}{2}}dt.\eqno(2.10)$$ \noindent
Thus
$$L'(\alpha_s)=-\int_{t=a}^bh(t,s)^{-\frac{1}{2}}A(t,s)dt;\eqno(2.11)$$
$$L''(\alpha_s)=-\int_{t=a}^b[h(t,s)^{-\frac{3}{2}}A(t,s)^2+h(t,s)^{-\frac{1}{2}}\frac{\partial
A(t,s)}{\partial s}]dt.\eqno(2.12)$$ \noindent So
$$L''(0)=-\int_{t=a}^b[A(t,0)^2+\frac{\partial
A(t,s)}{\partial s}(t,0)]dt.\eqno(2.13)$$ \noindent If $\gamma(t)$
is a unit timelike geodesic in $(M,g)$, then $A(t,0)=g(V',\gamma')$.
It remains to calculate $\frac{\partial A(t,s)}{\partial s}(t,0)$.
By (2.8) and commuting the differentiation, we have
$$\frac{\partial A(t,s)}{\partial s}(t,0)=g_B(\nabla^B_{\frac{\partial
\alpha_B}{\partial s }}\nabla^B_{\frac{\partial \alpha_B}{\partial
t}}\frac{\partial \alpha_B}{\partial s },\frac{\partial
\alpha_B}{\partial t })|_{s=0}+g_B(V_B',V_B')$$
$$+\frac{1}{2}\sum
_{i=1}^m\frac{\partial^2 }{\partial
s^2}(b_i^2(\alpha_B,\alpha_i))(t,s)|_{s=0}g_{F_i}(\gamma_{F_i}'(t),\gamma_{F_i}'(t))$$
$$
+2\sum _{i=1}^m\frac{\partial }{\partial
s}(b_i^2(\alpha_B,\alpha_i))(t,s)|_{s=0}g_{F_i}(V_{F_i}'(t),\gamma_{F_i}'(t))
$$
$$+\sum _{i=1}^mb_i^2(\gamma_{B},\gamma_{F_i})(t)
g_{F_i}(\nabla^{F_i}_{\frac{\partial \alpha_{F_i}}{\partial s
}}\nabla^{F_i}_{\frac{\partial \alpha_{F_i}}{\partial
t}}{\frac{\partial \alpha_{F_i}}{\partial s }},{\frac{\partial
\alpha_{F_i}}{\partial t }})|_{s=0} + \sum
_{i=1}^mb_i^2(\gamma_{B},\gamma_{F_i})(t)
g_{F_i}(V_{F_i}',V_{F_i}')$$

$$=g_B(V_B',V_B')-g_B(R^B(V_B,\gamma_B')\gamma_B',V_B)+g_B(\nabla^B_{\frac{\partial
\alpha_B}{\partial t }}\nabla^B_{\frac{\partial \alpha_B}{\partial
s}}\frac{\partial \alpha_B}{\partial s },\frac{\partial
\alpha_B}{\partial t })|_{s=0}$$
$$ +\frac{1}{2}\sum
_{i=1}^m\frac{\partial^2 }{\partial
s^2}(b_i^2(\alpha_B,\alpha_i))(t,s)|_{s=0}g_{F_i}(\gamma_{F_i}'(t),\gamma_{F_i}'(t))$$
$$+2\sum _{i=1}^m\frac{\partial }{\partial
s}(b_i^2(\alpha_B,\alpha_i))(t,s)|_{s=0}g_{F_i}(V_{F_i}'(t),\gamma_{F_i}'(t))
$$
$$-\sum _{i=1}^mb_i^2(\gamma_{B},\gamma_{F_i})(t)
g_{F_i}(R^{F_i}(V_{F_i},\gamma_{F_i}')\gamma_{F_i}',V_{F_i})$$

$$+\sum _{i=1}^mb_i^2(\gamma_{B},\gamma_{F_i})(t)
g_{F_i}(\nabla^{F_i}_{\frac{\partial \alpha_{F_i}}{\partial t
}}\nabla^{F_i}_{\frac{\partial \alpha_{F_i}}{\partial
s}}{\frac{\partial \alpha_{F_i}}{\partial s }},{\frac{\partial
\alpha_{F_i}}{\partial t }})|_{s=0}  + \sum
_{i=1}^mb_i^2(\gamma_{B},\gamma_{F_i})(t)
g_{F_i}(V_{F_i}',V_{F_i}').\eqno(2.14)$$ \noindent By the
proposition 2.10 and preserving the metric, we have
$$g_B(\nabla^B_{\frac{\partial
\alpha_B}{\partial t }}\nabla^B_{\frac{\partial \alpha_B}{\partial
s}}\frac{\partial \alpha_B}{\partial s },\frac{\partial
\alpha_B}{\partial t
})|_{s=0}=\frac{d}{dt}[g_B(\nabla^B_{V_B}V_B,\gamma_B')]-\frac{1}{2}\sum
_{i=1}^mg_B(\nabla^B_{V_B}V_B,{\rm
grad}_Bb_i^2)g_{F_i}(\gamma_{F_i}',\gamma_{F_i}'),\eqno(2.15)$$ We
note that
$$\frac{\partial^2 }{\partial
s^2}(b_i^2(\alpha_B,\alpha_i))(t,s)=(V_B^2+2V_BV_{F_i}+V_{F_i}^2)(b_i^2),\eqno(2.16)$$
\noindent so

$$g_B(\nabla^B_{\frac{\partial \alpha_B}{\partial t
}}\nabla^B_{\frac{\partial \alpha_B}{\partial s}}\frac{\partial
\alpha_B}{\partial s },\frac{\partial \alpha_B}{\partial t })|_{s=0}
 +\frac{1}{2}\sum _{i=1}^m\frac{\partial^2 }{\partial
s^2}(b_i^2(\alpha_B,\alpha_i))(t,s)|_{s=0}g_{F_i}(\gamma_{F_i}'(t),\gamma_{F_i}'(t))$$
$$=\frac{d}{dt}[g_B(\nabla^B_{V_B}V_B,\gamma_B')]+\frac{1}{2}\sum
_{i=1}^m{\rm
Hess}^B(b_i^2)(V_B,V_B)g_{F_i}(\gamma_{F_i}',\gamma_{F_i}')$$
$$+\frac{1}{2} \sum
_{i=1}^m(2V_BV_{F_i}+V_{F_i}^2)(b_i^2)g_{F_i}(\gamma_{F_i}',\gamma_{F_i}').\eqno(2.17)$$
\noindent By the proposition 2.10, similarly we can get\\
$$\sum _{i=1}^mb_i^2(\gamma_{B},\gamma_{F_i})(t)
g_{F_i}(\nabla^{F_i}_{\frac{\partial \alpha_{F_i}}{\partial t
}}\nabla^{F_i}_{\frac{\partial \alpha_{F_i}}{\partial
s}}{\frac{\partial \alpha_{F_i}}{\partial s }},{\frac{\partial
\alpha_{F_i}}{\partial t }})|_{s=0}$$ $$ =\frac{d}{dt}[\sum
_{i=1}^mb_i^2(\gamma_{B},\gamma_{F_i})(t)g_{F_i}(\nabla^{F_i}_{V_{F_i}}V_{F_i},\gamma_{F_i}')]-\frac{1}{2}
\sum _{i=1}^mg_{F_i}(\nabla^{F_i}_{V_{F_i}}V_{F_i},{\rm
grad}_{F_i}(b_i^2))g_{F_i}(\gamma_{F_i}',\gamma_{F_i}').\eqno(2.18)$$
\noindent By (2.13) (2.14),(2.17) and (2.18), we get\\

\noindent {\bf Proposition 2.12} {\it Let
$M=B\times_{b_1}F_1\times_{b_2}F_2\cdots
    \times_{b_m}F_m$ be a multiply twisted product. Let
    $\gamma:[a,b]\rightarrow (M,g)$ be a unit timelike geodesic and
    $\alpha:[a,b]\times (-\varepsilon,\varepsilon)\rightarrow (M,g)$
    be a smooth variation of $\gamma(t)$ with variation vector field
    $V=(V_B,V_{F_1},\cdots,V_{F_m}).$ Then}
    $$L''(0)=-\int_{t=a}^b\left\{g(V',\gamma')^2+g_B(V_B',V_B')-
g_B(R^B(V_B,\gamma_B')\gamma_B',V_B)\right.$$ $$+ \sum
_{i=1}^mb_i^2(\gamma_{B},\gamma_{F_i})(t)
[g_{F_i}(V_{F_i}',V_{F_i}')-g_{F_i}(R^{F_i}(V_{F_i},\gamma_{F_i}')\gamma_{F_i}',V_{F_i})]$$
$$
+\frac{1}{2}\sum _{i=1}^m[{\rm Hess}^B(b_i^2)(V_B,V_B)+{\rm
Hess}^{F_i}(b_i^2)(V_{F_i},V_{F_i})
+2V_BV_{F_i}(b_i^2)]g_{F_i}(\gamma_{F_i}',\gamma_{F_i}')$$ $$
\left.+ 2 \sum
_{i=1}^m(V_B+V_{F_i})(b_i^2)g_{F_i}(V_{F_i}',\gamma_{F_i}')\right\}dt$$
$$
-\left[g_B(\nabla^B_{V_B}V_B,\gamma_B')+ \sum
_{i=1}^mb_i^2(\gamma_{B},\gamma_{F_i})(t)g_{F_i}(\nabla^{F_i}_{V_{F_i}}V_{F_i},\gamma_{F_i}')\right]|_a^b.\eqno(2.19)$$\\

\indent In studying the second variation and index form, it suffices
to consider vector fields perpendicular to the given geodesic
$\gamma(t)$. Let $V^\bot(\gamma)$ denote the vector space of
piecewise smooth vector fields $V$ along $\gamma$ with
$g(V,\gamma')=0$ and let $V_0^\bot(\gamma)=\{V\in
V^\bot(\gamma)|V(a)=V(b)=0\}.$ Then guided by the result of
Proposition 2.12, the index form
$$I:~V_0^\bot(\gamma)\times V_0^\bot(\gamma)\rightarrow {\bf R}$$
\noindent should be given by\\
$$I(V,V)==-\int_{t=a}^b\left\{g_B(V_B',V_B')-
g_B(R^B(V_B,\gamma_B')\gamma_B',V_B)\right.$$ $$+ \sum
_{i=1}^mb_i^2(\gamma_{B},\gamma_{F_i})(t)
[g_{F_i}(V_{F_i}',V_{F_i}')-g_{F_i}(R^{F_i}(V_{F_i},\gamma_{F_i}')\gamma_{F_i}',V_{F_i})]$$
$$
+\frac{1}{2}\sum _{i=1}^m[{\rm Hess}^B(b_i^2)(V_B,V_B)+{\rm
Hess}^{F_i}(b_i^2)(V_{F_i},V_{F_i})
+2V_BV_{F_i}(b_i^2)]g_{F_i}(\gamma_{F_i}',\gamma_{F_i}')$$ $$
\left.+ 2 \sum
_{i=1}^m(V_B+V_{F_i})(b_i^2)g_{F_i}(V_{F_i}',\gamma_{F_i}')\right\}dt\eqno(2.20)$$
\noindent $I(V,W)$ could be obtained from (2.20) by polarization.\\
\indent Now we specialize to the index form to the case that
$(M,g)=({\bf R}\times {\bf R}\times F,
-du^2+f_1^2(u,x)dx^2+f_2^2(u)g_F)$ where $(F,g_F)$ is a Riemannian
manifold. We call such Lorentzian manifolds as static multiply
twisted product spacetimes. Let
$\gamma(t)=(\gamma_B(t),\gamma_1(t),\gamma_2(t))$ denote a unit
timelike geodesic segment. Consider variation vector fields
$V=(V_B,V_1,V_2)$ along $\gamma$ with $g(V,\gamma')=0$. We begin
with a special case in which the timelike geodesic $\gamma(t)$ is of
the form $\gamma(t)=(t,x_0,q)$. By  $g(V,\gamma')=0$,
then $V_B=0$. By $\gamma'_{F_i}=0$, we get\\

\noindent {\bf Proposition 2.13} {\it Let $(M,g)=({\bf R}\times {\bf
R}\times F, -du^2+f_1^2(u,x)dx^2+f_2^2(u)g_F)$ . Let
    $\gamma:[a,b]\rightarrow (M,g)$ be a unit timelike geodesic having form $\gamma(t)=(t,x_0,q)$ Then}
$I:~V_0^\bot(\gamma)\times V_0^\bot(\gamma)\rightarrow {\bf R}$ is
given by \\
$$I(V,V)=-\int_{t=a}^b\left[f_1^2(t,x_0)g_{R}(V_1',V_1')+f_2^2(t)g_F(V_2',V_2')\right].\eqno(2.21)$$

\indent By Proposition 2.10, we have\\

\noindent {\bf Corollary 2.14} {\it Let $(M,g)=({\bf R}\times {\bf
R}\times F, -du^2+f_1^2(u,x)dx^2+f_2^2(u)g_F)$ . Then
    $\gamma(t)=(\tau(t),\tau_1(t),\gamma_F(t))$ be a geodesic if and
    only if }\\
     $$(1)~\tau''(t)-\frac{1}{2}\tau_1'(t)^2\frac{\partial
    f_1^2}{\partial u}(\tau(t),\tau_1(t))-\frac{1}{2}\frac{\partial
    f_2^2}{\partial u}(\tau(t))
g_F(\gamma_F',\gamma_F')=0;\eqno(2.22)$$
$$(2)~\tau_1''(t)+\frac{1}{f_1^2(\tau(t),\tau_1(t))}\left[\frac{\partial
    f_1^2}{\partial
    u}(\tau(t),\tau_1(t))\tau'(t)\tau_1'(t)+\frac{1}{2}\frac{\partial
    f_1^2}{\partial
    x}(\tau(t),\tau_1(t))\tau_1'(t)^2\right]=0;\eqno(2.23)$$
$$(3)~\gamma_F''+\frac{1}{f_2^2(\tau(t))}\frac{\partial
    f_2^2}{\partial
    u}(\tau(t))\tau'(t)\gamma_{F}'(t)=0.\eqno(2.24)$$\\

    \indent Let $V_B=\nu_B(t)\frac{d}{du}|_{\tau(t)}$ and
    $V_1=\nu_1(t)\frac{d}{dx}|_{\tau_1(t)}$, then direct computation
    show that\\

\noindent {\bf Proposition 2.15} {\it Let $(M,g)=({\bf R}\times {\bf
R}\times F, -du^2+f_1^2(u,x)dx^2+f_2^2(u)g_F)$ . Let
    $\gamma:[a,b]\rightarrow (M,g)$ be a unit timelike geodesic
 and $V\in V_0^\bot(\gamma)$, then}\\
$$I(V,V)=-\int_{t=a}^b\left\{-\nu_B'(t)^2+f_1^2(\tau(t),\tau_1(t))\nu_1'(t)^2
+\frac{1}{2}\left[ \nu_B^2(t)\frac{\partial^2f_1^2}{\partial
u^2}(\tau(t),\tau_1(t))\right.\right.$$ $$\left.+
\nu_1^2(t)\frac{\partial^2f_1^2}{\partial x^2}(\tau(t),\tau_1(t))+
2\nu_B(t)\nu_1(t)\frac{\partial^2f_1^2}{\partial u
\partial x}(\tau(t),\tau_1(t))\right]\tau_1'(t)^2$$
$$+ 2\left[ \nu_B(t)\frac{\partial f_1^2}{\partial
u}(\tau(t),\tau_1(t))+\nu_1(t)\frac{\partial f_1^2}{\partial
x}(\tau(t),\tau_1(t))\right]\nu_1'(t)\tau_1'(t)$$
$$+f_2^2(\tau(t))\left[g_{F}(V_{2}',V_{2}')-g_{F}(R^{F}(V_{2},\gamma_{2}')\gamma_{2}',V_{2})\right]$$
$$\left.+\frac{1}{2}\nu_B^2(t)\frac{\partial^2f_2^2}{\partial
u^2}(\tau(t))g_F(\gamma_2'(t),\gamma_2'(t))+2 \nu_B(t)\frac{\partial
f_2^2}{\partial
u}(\tau(t))g_F(V_2'(t),\gamma_2'(t))\right\}dt.\eqno(2.25)$$

\section{Killing vector fields
 on the multiply twisted product manifolds}

\quad In this section, we develop the propertied of Killing vector
fields, then we focus our attention on some special cases and
characterize Killing vector fields on these spaces.\\

\noindent {\bf Lemma 3.1} {\it  Let
$M=B\times_{b_1}F_1\times_{b_2}F_2\cdots
    \times_{b_m}F_m$ be a multiply twisted product. Let $X\in
    \Gamma(TB),~U\in \Gamma(TF_i)$, then}\\
$$L^M_Xg_M=L_X^Bg_B+\sum_{i=1}^mX(b_i^2)g_{F_i};~~L^M_Ug_M=b_i^2L_U^{F_i}g_{F_i}+U(b_i^2)g_{F_i}.\eqno(3.1)$$\\

\noindent {\bf Proposition 3.2} {\it  Let
$M=B\times_{b_1}F_1\times_{b_2}F_2\cdots
    \times_{b_m}F_m$ be a multiply twisted product. Let $X\in
    \Gamma(TB),~U_i\in \Gamma(TF_i)$, then $K=X+U_1+\cdots U_m$ is a Killing vector fields if and only if
    $X$ is a Killing vector fields on $B$ and $U_i$ is a conformal Killing vector
    fields on $F_i$ with the conformal factor
    $-\frac{X(b_i^2)+U_i(b_i^2)}{b_i^2}$.}\\

\noindent {\bf Lemma 3.3} {\it  Let
$M=B\times_{b_1}F_1\times_{b_2}F_2\cdots
    \times_{b_m}F_m$ be a multiply twisted product. If $K$ is a
    Killing vector field on $M$, then $K_B=K(\cdot,q_1,\cdots,q_m)$
    is a Killing vector field on $B$, and
    $K_{F_i}=K(p,q_1,\cdots,q_{i-1},\cdot,q_{i+1},\cdots,q_m)$ is Killing vector field on
    $F_i$.}\\

Nextly, we assume $m=2$ and $b_i=f_i$. Let $K_a,~a\in \{1,\cdots
,\overline{m_B}\}$ be a basis of Killing vector fields of $(B,g_B)$
and $W_b,~b\in\{1,\cdots, \overline{m_{F_1}}\}$ be a basis of
conformal Killing vector fields of $(F_1,g_{F_1})$ with
$L_{W_b}^{F_1}g_{F_1}=2\sigma_bg_{F_1}$ and $G_c,~c\in\{1,\cdots,
\overline{m_{F_2}}\}$ be a basis of conformal Killing vector fields
of $(F_2,g_{F_2})$ with
$L_{G_c}^{F_2}g_{F_2}=2\overline{\sigma_c}g_{F_2}$. Let $x,y,z$
denote variables on $B,F_1,F_2$ respectively. If $K$ is a Killing
vector field on $M$, then by Lemma 3.3, there are functions
$\mu^a(y,z),~\delta^b(x,z),\lambda^c(x,y)$ such that
$$K=K_B+K_{F_1}+K_{F_2},~K_B=\mu^aK_a,~K_{F_1}=\delta^bW_b,~K_{F_2}=\lambda^cG_c.\eqno(3.2)$$
We note that if $\phi,Z$ and $T$ are respectively a function, a
vector field and a 2-covariant tensor field on $M$, then\\
$$L_{\phi Z}T(\cdot,\cdot)=\phi
L_ZT(\cdot,\cdot)+d\phi(\cdot)\otimes T(Z,\cdot)+T(Z,\cdot)\otimes
d\phi(\cdot).\eqno(3.3)$$ \noindent Let
$\hat{K_a}=g_B(K_a,\cdot),~\hat{W_b}=g_{F_1}(W_b,\cdot)~\hat{G_c}=g_{F_2}(G_c,\cdot).$
 By Lemma 3.1 and (3.2),
(3.3), we get\\
$$L_Kg_M=2f_1^2\left[K_B(lnf_1)+K_{F_1}(lnf_1)+\delta^b\sigma_b\right]g_{F_1}+2f_2^2\left[K_B(lnf_2)+K_{F_2}(lnf_2)+\lambda^c\overline{\sigma_c}
\right]g_{F_2}$$ $$+d\mu^a\otimes \widehat{K_a}+\widehat{K_a}\otimes
d\mu^a+f_1^2(d\delta^b\otimes \widehat{W_b}+\widehat{W_b}\otimes
d\delta^b)+f_2^2(d\lambda^c\otimes
\widehat{G_c}+\widehat{G_c}\otimes d\lambda^c)$$
$$=2f_1^2\left[K_B(lnf_1)+K_{F_1}(lnf_1)+\delta^b\sigma_b\right]g_{F_1}+2f_2^2\left[K_B(lnf_2)+K_{F_2}(lnf_2)+\lambda^c\overline{\sigma_c}
\right]g_{F_2}$$ $$+ (d_{F_1}\mu^a\otimes
\widehat{K_a}+f_1^2\widehat{W_b}\otimes d_B\delta^b)+(\widehat{K_a}
\otimes d_{F_1}\mu^a+f_1^2d_B\delta^b\otimes \widehat{W_b})$$ $$ +
(d_{F_2}\mu^a\otimes \widehat{K_a}+f_2^2\widehat{G_c}\otimes
d_B\lambda^c)+(\widehat{K_a}\otimes
d_{F_2}\mu^a+f_2^2d_B\lambda^c\otimes \widehat{G_c})$$
$$+(f_1^2d_{F_2}\delta^b\otimes
\widehat{W_b}+f_2^2\widehat{G_c}\otimes
d_{F_1}\lambda^c)+(f_1^2\widehat{W_b}\otimes
d_{F_2}\delta^b+f_2^2d_{F_1}\lambda^c\otimes\widehat{G_c}).\eqno(3.4)$$
So we obtain the following result.\\

\noindent {\bf Proposition 3.4} {\it  Let
$M=B\times_{f_1}F_1\times_{f_2}F_2$. Let $K$ be a vector field on
$M$ as in (3.2), then $K$ is a Killing vector field if and only if
the following relations hold:}\\
$$K_B(lnf_1)+K_{F_1}(lnf_1)+\delta^b\sigma_b=0;~K_B(lnf_2)+K_{F_2}(lnf_2)+\lambda^c\overline{\sigma_c}=0;\eqno(3.5)$$
$$d_{F_1}\mu^a\otimes
\widehat{K_a}+f_1^2\widehat{W_b}\otimes
d_B\delta^b=0;~d_{F_2}\mu^a\otimes
\widehat{K_a}+f_2^2\widehat{G_c}\otimes d_B\lambda^c=0;\eqno(3.6)$$
$$f_1^2d_{F_2}\delta^b\otimes
\widehat{W_b}+f_2^2\widehat{G_c}\otimes
d_{F_1}\lambda^c=0.\eqno(3.7)$$\\

\indent Now we let $M=I\times_{f_1}F_1\times_{f_2}F_2$ with the
metric tensor $-dt^2+f_1(t)^2g_{F_1}+f_2(t)^2g_{F_2}.$ By (3.2),
then $K=\mu (y,z)\frac{\partial}{\partial
t}+\delta^b(t,z)W_b+\lambda^c(t,y)G_c.$ By Proposition 3.4, we have
if $K$ is a Killing vector field, then \\
$$K_B(lnf_1)+\delta^b\sigma_b=0;~K_B(lnf_2)+\lambda^c\overline{\sigma_c}=0;\eqno(3.8)$$
$$-d_{F_1}\mu \otimes
dt+f_1^2\widehat{W_b}\otimes d_B\delta^b=0;~-d_{F_2}\mu\otimes dt
+f_2^2\widehat{G_c}\otimes d_B\lambda^c=0;\eqno(3.9)$$
$$f_1^2d_{F_2}\delta^b\otimes
\widehat{W_b}+f_2^2\widehat{G_c}\otimes
d_{F_1}\lambda^c=0.\eqno(3.10)$$ \noindent By (3.9), then
$$d_{F_1}\mu(y,z)=f_1^2\frac{\partial \delta^b(t,z)}{\partial
t}\widehat{W_b}.\eqno(3.11)$$ Since $\widehat{W_b}$ is a base, by
separation of variables, we obtain\\
$$\delta^b(t,z)=\psi_b(z)\int_{t=t_0}^tf_1^{-2}(u)du+\widetilde{\psi_b(z)},~~
d_{F_1}\mu(y,z)=\psi_b(z)\widehat{W_b}.\eqno(3.12)$$ Similarly, we
obtain\\
$$d_{F_2}\mu(y,z)=f_2^2\frac{\partial \lambda^c(t,y)}{\partial
t}\widehat{G_c}=\phi_c(y)\widehat{G_c};~~\lambda^c(t,y)=\phi_c(y)\int_{t=t_0}^tf_2^{-2}(u)du+\widetilde{\phi_c(y)}.\eqno(3.13)$$
One derives (3.8) with respect to $t$, then
$$\mu(y,z)\left(\frac{f_1'}{f_1}\right)'+f_1^{-2}(t)\psi_b(z)\sigma_b(y)=0\eqno(3.14)$$
\noindent {\bf Case I)} $\mu=0$, by (3.12) and (3.13),
$$\psi_b(z)=0,~\delta^b(t,z)=\widetilde{\psi_b(z)},~\phi_c(y)=0,~
\lambda^c(t,y)=\widetilde{\phi_c(y)}
,\widetilde{\psi_b(z)}\sigma_b=0,~\widetilde{\phi_c(y)}\overline{\sigma_c}=0.\eqno(3.15)$$
So we have\\

\noindent {\bf Proposition 3.5} {\it Let
$M=I\times_{f_1}F_1\times_{f_2}F_2$ with the metric tensor
$-dt^2+f_1(t)^2g_{F_1}+f_2(t)^2g_{F_2}.$ $K$ is a Killing vector
field with $\mu=0$, then $K=K_{F_1}(y,z)+K_{F_2}(y,z)$ and
$K_{F_1}(y,z_0)$ is a Killing vector field on $(F_1,g_{F_1})$ and
$K_{F_2}(y_0,z)$ is a Killing vector field on $(F_2,g_{F_2}).$}\\

\noindent \noindent {\bf Case II)} $\mu\neq 0$. By (3.14), then
$$\left(\frac{f_1'}{f_1}\right)'f_1^{2}(t)=-\frac{\psi_b(z)\sigma_b(y)}{\mu(y,z)}=C_F^1;\eqno(3.16)$$
$$\left(\frac{f_2'}{f_2}\right)'f_2^{2}(t)=-\frac{\phi_c(y)\overline{\sigma_c(z)}}{\mu(y,z)}=C_F^2,\eqno(3.17)$$
\noindent where $C_F^1,C_F^2$ are constants. By (3.2),(3.12),(3.13),(3.16) and (3.17), we obtain\\
$$K=\mu(y,z)\frac{\partial}{\partial t}+{\rm
grad}_{F_1}\mu\int_{t_0}^tf_1^{-2}(u)du+{\rm
grad}_{F_2}\mu\int_{t_0}^tf_2^{-2}(u)du+\widetilde{\psi_b}(z)W_b+\widetilde{\phi_c}(y)G_c;\eqno(3.18)$$
$${\rm Hess}^\mu_{F_1}+C_F^1\mu g_{F_1}=0;~{\rm Hess}^\mu_{F_2}+C_F^2\mu g_{F_2}=0;\eqno(3.19)$$
We may assume that $\mu$ is not a constant. One can derive (3.8)
with respect to $t$ and by (3.12) and (3.16), we get
$$-\psi_b\sigma_b(lnf_1)'(t_0)+C_F^1\widetilde{\psi_b}(z)\sigma_b(y)=0.\eqno(3.20)$$
So if ${\bf C_F^1\neq 0}$, then
$T^{F_1}=\widetilde{\psi_b}(z)W_b-\frac{(lnf_1)'(t_0)}{C_F^1}{\psi_b}(z)W_b$
is a Killing vector field on $F_1$ for fixed $z$ and
$$\widetilde{\psi_b}(z)W_b=T^{F_1}+\frac{(lnf_1)'(t_0)}{C_F^1}{\rm
grad}_{F_1}\mu;\eqno(3.21)$$ \noindent Similarly, if ${\bf C_F^2\neq
0}$, then
$$\widetilde{\phi_c}(y)G_c=T^{F_2}+\frac{(lnf_2)'(t_0)}{C_F^2}{\rm
grad}_{F_2}\mu.\eqno(3.22)$$ By (3.16),(3.17),(3.18), (3.21) and
(3.22), we
obtain\\

\noindent {\bf Proposition 3.7} {\it If
$M=I\times_{f_1}F_1\times_{f_2}F_2$ with the metric tensor
$-dt^2+f_1(t)^2g_{F_1}+f_2(t)^2g_{F_2}$ admits a Killing vector
field with $\mu\neq 0$, then
$\left(\frac{f_i'}{f_i}\right)'f_i^{2}(t)=C_F^i~for ~i=1,2$.}\\

\noindent {\bf Proposition 3.8} {\it Let
$M=I\times_{f_1}F_1\times_{f_2}F_2$ with the metric tensor
$-dt^2+f_1(t)^2g_{F_1}+f_2(t)^2g_{F_2}.$ and
$\left(\frac{f_i'}{f_i}\right)'f_i^{2}(t)=C_F^i\neq 0~for ~i=1,2$.
Then its Killing vector field is given by }\\
$$K=\mu(y,z)\frac{\partial}{\partial t}+{\rm
grad}_{F_1}\mu\left(\int_{t_0}^tf_1^{-2}(u)du+
\frac{(lnf_1)'(t_0)}{C_F^1}\right)$$ $$ +{\rm
grad}_{F_2}\mu\left(\int_{t_0}^tf_2^{-2}(u)du+\frac{(lnf_2)'(t_0)}{C_F^2}\right)
+T^{F_1}(y,z)+T^{F_2}(y,z),\eqno(3.23)$$
 \noindent {\it where
$T^{F_1}(y,z_0)$ is a Killing vector field on $(F_1,g_{F_1})$ and
$T^{F_2}(y_0,z)$ is a Killing vector field on $(F_2,g_{F_2})$ and $\mu$ satisfies (3.19).}\\

\indent When $C_F^i=0$, we may obtain similar results like Theorem
4.7 in [Sa].\\
\indent Nextly, we consider {\bf the Killing vector fields on
${\mathbb{R}}^2\times F$ with the metric tensor
$-dt^2+ds^2+f^2(t)g_F$.}\\
 \indent Direct computations show that
vector fields $K_1=\frac{\partial}{\partial
t},~K_2=\frac{\partial}{\partial s},~K_3=s\frac{\partial}{\partial
t}+t\frac{\partial}{\partial s}$ are the basis of Killing vector
fields on $\textrm{R}^2$ with the metric tensor $-dt^2+ds^2$. Let
$W_b,~b\in\{1,\cdots, \overline{m_{F}}\}$ be a basis of conformal
Killing vector fields of $(F,g_{F})$ with
$L_{W_b}^{F}g_{F}=2\sigma_bg_{F}$, then
$$K=\mu^aK_a+\delta^bW_b;~
K_B(lnf)+\delta^b\sigma_b=0;~d\mu^a\otimes
\widehat{K_a}+f^2\widehat{W_b}\otimes d\delta^b=0.\eqno(3.24)$$ By
$\widehat{K_1}=-dt,~\widehat{K_2}=ds,~\widehat{K_3}=-sdt+tds$ and
(3.24), we have:
$$d\mu^1+sd\mu^3-f^2(t)\frac{\partial \delta^b}{\partial
t}\widehat{W_b}=0;\eqno(3.25)$$
$$~d\mu^2+td\mu^3+f^2(t)\frac{\partial \delta^b}{\partial
s}\widehat{W_b}=0.\eqno(3.26)$$ We
 derive (3.25) with respect to $s$, then
$$\delta^b=C_bs\int_{t_0}^tf^{-2}(u)du+\widehat{C_b}(s)+L_b(t),\eqno(3.27)$$
\noindent where $C_b$ is a constant. By (3.25) and (3.27), we have
$$d\mu^3=C_b\widehat{W_b},~d\mu^1=d_b\widehat{W_b},~\delta^b=(C_bs+d_b)\int_{t_0}^tf^{-2}(u)du+\widetilde{C_b}(s),\eqno(3.28)$$
\noindent where $d_b$ is a constant. By (3.26) and (3.28), we obtain
$$d\mu^2=\left(-C_bt-f^2(t)C_b\int_{t_0}^tf^{-2}(u)du-f^2(t)\widetilde{C_b}'(s)\right)\widehat{W_b}.\eqno(3.29)$$
So
$$d\mu^2=e_b\widehat{W_b},~C_bt+f^2(t)C_b\int_{t_0}^tf^{-2}(u)du+f^2(t)\widetilde{C_b}'(s)=e_b\eqno(3.30)$$
\noindent where $e_b$ is a constant. Then
$$\widetilde{C_b}'(s)=-C_btf^{-2}(t)-C_b\int_{t_0}^tf^{-2}(u)du+f^{-2}(t)e_b,\eqno(3.31)$$
and
$$\widetilde{C_b}'(s)=\widehat{e_b},~~\widetilde{C_b}(s)=\widehat{e_b}s+\widetilde{e_b},\eqno(3.32)$$
\noindent where $=\widehat{e_b}, \widetilde{e_b}$ are constants. By
(3.28) and (3.32), we have
$$\delta^b=(C_bs+d_b)\int_{t_0}^tf^{-2}(u)du+\widehat{e_b}s+\widetilde{e_b}.\eqno(3.33)$$
We derive (3.31) with respect to $t$, then
$$C_bf=f'(C_bt-e_b).\eqno(3.34)$$
\indent {\bf Case I)} There is a $C_b\neq 0$ i.e. $\mu^3$ is not a
constant. So by (3.34), then $f=A(t-\frac{e_b}{C_b})=At+B~for ~A\neq
0$ and $f$ is fixed and $\frac{e_b}{C_b}$ is a constant independent
of $b$. So by (3.28) and (3.30),
$$d\mu^2=\lambda d\mu^3,~~\mu^2=\lambda
\mu^3+\widetilde{\lambda}.\eqno(3.35)$$ \noindent where $\lambda,
\widetilde{\lambda}$ are constants. By the second equation in (3.24)
and (3.33), we can obtain
$$\frac{A}{At+B}\mu_1+\left[d_b\int_{t_0}^tf^{-2}(u)du+\widetilde{e_b}\right]\sigma_b=0,\eqno(3.36)$$
$$\frac{A}{At+B}\mu_3+\left[C_b\int_{t_0}^tf^{-2}(u)du+\widehat{e_b}\right]\sigma_b=0.\eqno(3.37)$$
\noindent Derive (3.36) and (3.37) with respect to $t$, then
$$\mu_1-\frac{1}{A^2}d_b\sigma_b=0;~~\mu_3-\frac{1}{A^2}C_b\sigma_b=0.\eqno(3.38)$$
So
$${\rm Hess}_F^{\mu_1}=A^2\mu_1g_F;~~{\rm
Hess}_F^{\mu_3}=A^2\mu_3g_F.\eqno(3.39)$$ By (3.31),(3.32) and
$f=At+B$, then $$\widehat{e_b}=-\frac{C_b}{A(At_0+B)}.\eqno(3.40)$$
By (3.24),(3.33),(3.35) and (3.39), we obtain
$$K=(\mu_1+s\mu_3)\frac{\partial}{\partial
t}+[(\lambda+t)\mu_3+\widetilde{\lambda}]\frac{\partial}{\partial
s}-\frac{s}{A(At+B)}{\rm grad}_F\mu_3$$
$$+\left[-\frac{1}{A(At+B)}+\frac{1}{A(At_0+B)}\right]{\rm
grad}_F\mu_1+W_*,\eqno(3.41)$$ \noindent where
$W_*=\widetilde{e_b}W_b$ is a conformal killing vector field on
$(F,g_F)$. In (3.36), we set $t=t_0$ and using (3.38), then
$$\frac{d_b\sigma_b}{A(At_0+B)}+\widetilde{e_b}\sigma_b=0.\eqno(3.42)$$
and $T=\frac{d_bW_b}{A(At_0+B)}+\widetilde{e_b}W_b$ is a Killing
vector field on $F$, So
$$K=(\mu_1+s\mu_3)\frac{\partial}{\partial
t}+[(\lambda+t)\mu_3+\widetilde{\lambda}]\frac{\partial}{\partial
s}-\frac{s}{A(At+B)}{\rm grad}_F\mu_3$$
$$-\frac{1}{A(At+B)}{\rm
grad}_F\mu_1+T,\eqno(3.43)$$
\\

\noindent {\bf Theorem 3.9}{ \it Let $M=\textrm{R}^2\times F$ with
the metric tensor $-dt^2+ds^2+f^2(t)g_F$ and $K$ is a Killing vector
field given by (3.24) and $\mu_3$ is not a constant, then
$f=At+B,~A\neq 0$ and $K$ can be expressed by (3.43) and
$\mu_1,\mu_3$ satisfy (3.39).}\\

\indent {\bf Case II)} $C_b=0$ for any $b$ and there is a $e_b\neq
0$ i.e. {\bf $\mu_3=k_0$ is a constant and  $\mu_2$ is not a
constant.} By (3.34), $f=l_0$ is a constant. By (3.33), then
$$\delta^b=\frac{d_b}{l_0^2}(t-t_0)+\widehat{e_b}s+\widetilde{e_b}.\eqno(3.44)$$
By (3.24) and (3.44), we obtain
$$d_b\sigma_b=0,~~\widetilde{e_b}\sigma_b=0,~\widehat{e_b}\sigma_b=0,\eqno(3.45)$$
so $\delta^bW_b=K_1t+K_2s+K_3$ where $K_1,K_2,K_3$ are Killing
vector fields on $F$. We also get by (3.31)
$${\rm Hess}^{\mu_1}_F={\rm Hess}^{\mu_2}_F=0.\eqno(3.46)$$\\

\noindent {\bf Theorem 3.10}{ \it Let $M=\textrm{R}^2\times F$ with
the metric tensor $-dt^2+ds^2+f^2(t)g_F$ and $K$ is a Killing vector
field given by (3.24) and $\mu_3=k_0$ and $\mu_2$ is not a constant,
then $f=l_0$ and $K$ can be expressed by }
$$K=(\mu_1+sk_0)\frac{\partial}{\partial t}+(\mu_2+tk_0)\frac{\partial}{\partial
s}+K_1t+K_2s+K_3,\eqno(3.47)$$ {\it where $K_1,K_2,K_3$ are Killing
vector fields on $F$ and
$\mu_1,\mu_3$ satisfy (3.46).}\\

\indent {\bf Case III)} $C_b=e_b=0$ for any $b$ i.e. {\bf
$\mu_3=k_0,~\mu_2=\overline{k_0}$ are constants}. In this case,
$$\delta^b={d_b}\int_{t_0}^tf^{-2}(u)du+\widetilde{e_b}.\eqno(3.48)$$
By (3.24), then
$$\mu_1\frac{f'}{f}+[{d_b}\int_{t_0}^tf^{-2}(u)du+\widetilde{e_b}]\sigma_b=0;~~k_0\frac{f'}{f}=0.\eqno(3.49)$$
When $k_0=0$, similar to the discussions in [Sa], we can obtain\\

\noindent {\bf Theorem 3.11}{ \it Let $M=\textrm{R}^2\times F$ with
the metric tensor $-dt^2+ds^2+f^2(t)g_F$ and $K$ is a Killing vector
field given by (3.24) and $\mu_3=0$ and $\mu_2=\overline{k_0}$, then
if $\mu_1=0$, $K$ can be expressed by }
$$K=\overline{k_0}\frac{\partial}{\partial s}+W_*,\eqno(3.50)$$
{\it where $W_*$ is a Killing vector field on $F$. If $\mu_1\neq 0$,
then $\left(\frac{f'}{f}\right)'f^2=C_F$ and when $C_F\neq 0$, then
$K$ can be expressed by}
$$K=\mu_1\frac{\partial}{\partial t}+\overline{k_0}\frac{\partial}{\partial
s}+{\rm
grad}_F\mu_1\left(\int_{t_0}^tf^{-2}(u)du+\frac{(lnf)'(t_0)}{C_F}\right)+T,\eqno(3.51)$$
{\it where $T$ is a Killing vector field on $F$ and ${\rm
Hess}_F^{\mu_1}+\mu_1C_Fg_F=0.$ If $C_F=0$ and $\mu_1$ is a nonzero
constant, then $K=\mu_1\frac{\partial}{\partial
t}+\overline{k_0}\frac{\partial}{\partial s}+W_*,$ where $W_*$ is
homothetic.}\\

\noindent {\bf Theorem 3.12}{ \it Let $M=\textrm{R}^2\times F$ with
the metric tensor $-dt^2+ds^2+f^2(t)g_F$ and $K$ is a Killing vector
field given by (3.24) and $\mu_3=k_0\neq 0$, $\mu_2=\overline{k_0}$,
then $f=l_0$ and $K$ can be expressed by }
$$K=(\mu_1+sk_0)\frac{\partial}{\partial t}+(\overline{k_0}+tk_0)\frac{\partial}{\partial
s}+\frac{{\rm grad}_F\mu_1}{l_0^2}(t-t_0)+T,\eqno(3.52)$$ {\it where
$T$ is a Killing vector field on $F$ and ${\rm
Hess}_F^{\mu_1}=0.$}\\

\noindent{\bf Remark.} Here we can not consider
$M=\textrm{R}^2\times F$ as $R\times (R\times _fF)$ with the metric
tensor $ds^2+(-dt^2+f^2(t)g_F)$. Although the space of conformal
Killing vector fields on $(F,g_F)$ is finite dimensional, but the
space of  Killing vector fields on $(R\times _fF,-dt^2+f^2(t)g_F)$
maybe is infinite dimensional.\\

 \indent Nextly we consider $M=I\times _f F$
with the metric tensor $-f_1^2(t)dt^2+f^2(t)g_F$ which generalizes
the generalized Robertson-Walker spacetime. We note that
$K_a=\frac{1}{f_1(t)}\frac{\partial}{\partial t}$ is the base of the
Killing vector fields on $(I,-f_1^2(t)dt^2)$ and
$\widehat{K_a}=-f_1(t)dt.$ Similar to the discussions in [Sa], we
get\\

\noindent {\bf Proposition 3.13}{ \it Let $M=I\times _f F$ with the
metric tensor $-f_1^2(t)dt^2+f^2(t)g_F$ and if $M$ admits a
non-trivial Killing vector field then
$\left(\frac{f'}{ff_1}\right)'\frac{f^2}{f_1}=C_F$. When $C_F\neq
0$, then $K$ can be expressed by}
$$K=\frac{\mu(x)}{f_1(t)}\frac{\partial}{\partial
t}+\left(\frac{(lnf)'(t_0)}{C_Ff_1(t_0)}+\int_{t_0}^tf_1f^{-2}du\right){\rm
grad}_F\mu+T,\eqno(3.53)$$ {\it  where $T$ is a Killing vector field
on $F$ and ${\rm Hess}_F^{\mu}+\mu C_Fg_F=0.$ When $C_F=0$, if $\mu$
is a constant, then $K=\frac{\mu_0}{f_1(t)}\frac{\partial}{\partial
t}+W^*$ where $W^*$ is homothetic. Otherwise $\nabla^F\mu$ is a
non-zero parallel field.}\\

\indent We recall the definition of the curl operator on
semi-Riemannian manifolds, namely: if $V$ is a vector field on a
seni-Riemannian manifold $M$, then ${\rm curl}V$ is the
antisymmetric 2-covariant tensor defined by
$${\rm curl}V(X,Y):=g_M(\nabla_XV,Y)-g_M(\nabla_YV,X),\eqno(3.54)$$
where $X,Y\in \Gamma(TM)$. A vector field $V$ on a semi-Riemannian
manifold $M$ is said to be non-rotating if ${\rm curl}V(X,Y)=0$ for
all $X,Y\in \Gamma(TM)$. By the remark 5.1 in [DU], we know that $V$
is non-rotating iff it is parallel. By Proposition 2.2, then for
$X,Y\in \Gamma(TR)$ and $V,W\in\Gamma(TF)$\\
\noindent$(1) ~~\nabla_XY=\nabla^B_XY.$\\
      \noindent$(2)~~\nabla_XW=\nabla_WX=\frac{X(f)}{f}W.$\\
      \noindent$(3)~~\nabla_VW=\frac{ff'}{f_1^2}g_F(V,W)+\nabla^{F}_VW.$\\
We take a Killing vector field
$$K=\frac{\mu(x)}{f_1(t)}\frac{\partial}{\partial
t}+\int_{t_0}^tf_1f^{-2}du{\rm grad}_F\mu+W_*,\eqno(3.55)$$ Then $K$
is non-rotating iff $\nabla_{\frac{\partial}{\partial t}}K=0$ and
$\nabla_VK=0$ for all $V\in \Gamma(TF)$. Direct computations show
that
$$\nabla_{\frac{\partial}{\partial t}}\left[\frac{\mu(x)}{f_1(t)}\frac{\partial}{\partial
t}\right]=2\frac{f_1'(t)}{f_1^2}\mu(x)\frac{\partial}{\partial
t};~~~~\nabla_{\frac{\partial}{\partial
t}}W_*=\frac{f'}{f}W_*;\eqno(3.56)$$
$$\nabla_{\frac{\partial}{\partial t}}\left[\int_{t_0}^tf_1f^{-2}du{\rm
grad}_F\mu\right]=f_1f^{-2}{\rm
grad}_F\mu+\int_{t_0}^tf_1f^{-2}du\frac{f'}{f}{\rm
grad}_F\mu.\eqno(3.57)$$ By $\nabla_{\frac{\partial}{\partial
t}}K=0$ and (3.56), (3.57) and $\mu\neq 0$,  we obtain if $f'\neq
0$, $f_1=c_0$ is a constant and
$$\int_{t_0}^tf_1f^{-2}du{\rm grad}_F\mu+W_*=-\frac{f_1}{ff'}{\rm
grad}_F\mu,\eqno(3.58)$$ and
$$K=\frac{\mu(x)}{f_1(t)}\frac{\partial}{\partial
t}-\frac{f_1}{ff'}{\rm grad}_F\mu.\eqno(3.59)$$ So
$$0=\nabla_VK=\frac{\mu f'}{ff_1}V-\frac{f_1}{ff'}\nabla^F_V({\rm
grad}_F\mu).\eqno(3.60)$$ By ${\rm Hess}_F^{\mu}+\mu C_Fg_F=0$ and
(3.60), we have $(f')^2=-f_1^2C_F,~f=\pm f_1\sqrt{-C_F}t+r_0.$\\

\noindent {\bf Proposition 3.14}{ \it Let $M=I\times _f F$ with the
metric tensor $-f_1^2(t)dt^2+f^2(t)g_F$ and if $K$ is a non-trivial
non-rotating Killing vector field and $f'\neq 0$, then $f_1$ is a
constant and $f=\pm f_1\sqrt{-C_F}t+r_0$ and $K$ can be expressed
by}
$$K=\frac{\mu(x)}{f_1(t)}\frac{\partial}{\partial
t}+\frac{{\rm grad}_F\mu}{f_1C_Ft+r_0}.\eqno(3.61)$$

\section{Multiply twisted
 products with a semi-symmetric metric connection}

\noindent {\bf 4.1 Preliminaries}\\

\indent Let $M$ be a Riemannian manifold with Riemannian metric $g$.
A linear connection $\overline{\nabla}$ on a Riemannian manifold $M$
is called a {\it semi-symmetric connection} if the torsion tensor
$T$ of the connection $\overline{\nabla}$
$$T(X,Y)=\overline{\nabla}_XY-\overline{\nabla}_YX-[X,Y]\eqno(4.1)$$
satisfies
$$T(X,Y)=\pi(Y)X-\pi(X)Y,\eqno(4.2)$$
where $\pi$ is a 1-form associated with the vector field $P$ on $M$
defined by $\pi(X)=g(X.P).$ $\overline{\nabla}$ is called a {\it
semi-symmetric metric connection} if it satisfies
$\overline{\nabla}g=0.$ If $\nabla$ is the Levi-Civita connection of
$M$, the semi-symmetric metric connection $\overline{\nabla}$ is
given by
$$\overline{\nabla}_XY=\nabla_XY+\pi(Y)X-g(X,Y)P,\eqno(4.3)$$
(see [Ya]). Let $R$ and $\overline{R}$ be the curvature tensors of
$\nabla$ and $\overline{\nabla}$ respectively. Then $R$ and
$\overline{R}$ are related by
$$\overline{R}(X,Y)Z=R(X,Y)Z+g(Z,\nabla_XP)Y-g(Z,\nabla_YP)X$$
$$+g(X,Z)\nabla_YP-g(Y,Z)\nabla_XP+\pi(P)[g(X,Z)Y-g(Y,Z)X]$$
$$+[g(Y,Z)\pi(X)-g(X,Z)\pi(Y)]P+\pi(Z)[\pi(Y)X-\pi(X)Y],\eqno(4.4)$$
for any vector fields $X,Y,Z$ on $M$ [Ya]. By (4.3) and Proposition
2.2, we have\\

 \noindent {\bf Proposition 4.1} {\it Let $M=B\times_{b_1}F_1\times_{b_2}F_2\cdots
    \times_{b_m}F_m$ be a multiply twisted product and let $X,Y\in \Gamma(TB)$
      and $U\in \Gamma(TF_i)$, $W\in \Gamma(TF_j)$ and $P\in \Gamma(TB)$ . Then}\\
      \noindent$(1) ~~\overline{\nabla}_XY=\overline{\nabla}^B_XY.$\\
      \noindent$(2)~~\overline{\nabla}_XU=\frac{X(b_i)}{b_i}U.$\\
\noindent$(3)~~\overline{\nabla}_UX=[\frac{X(b_i)}{b_i}+\pi(X)]U.$\\
            \noindent$(4)~~\overline{\nabla}_UW=0~~if~ i\neq j.$\\
      \noindent$(5)~~\overline{\nabla}_UW=U(lnb_i)W+W(lnb_i)U-\frac{g_{F_i}(U,W)}{b_i}{\rm
      grad}_{F_i}b_i-b_ig_{F_i}(U,W){\rm
      grad}_{B}b_i+\nabla^{F_i}_UW-g(U,W)P~~if~ i= j.$\\

\noindent {\bf Proposition 4.2} {\it Let
$M=B\times_{b_1}F_1\times_{b_2}F_2\cdots
    \times_{b_m}F_m$ be a multiply twisted product and let $X,Y\in \Gamma(TB)$
      and $U\in \Gamma(TF_i)$, $W\in \Gamma(TF_j)$ and $P\in \Gamma(TF_k)$ . Then}\\
      \noindent$(1) ~~\overline{\nabla}_XY={\nabla}^B_XY-g(X,Y)P.$\\
      \noindent$(2)~~\overline{\nabla}_XU=\frac{X(b_i)}{b_i}U+g(P,U)X.$\\
\noindent$(3)~~\overline{\nabla}_UX=\frac{X(b_i)}{b_i}U.$\\
            \noindent$(4)~~\overline{\nabla}_UW=g(W,P)U~~if~ i\neq j.$\\
      \noindent$(5)~~\overline{\nabla}_UW=U(lnb_i)W+W(lnb_i)U-\frac{g_{F_i}(U,W)}{b_i}{\rm
      grad}_{F_i}b_i-b_ig_{F_i}(U,W){\rm
      grad}_{B}b_i+\nabla^{F_i}_UW+\pi(W)U-g(U,W)P~~if~ i= j.$\\

\indent By (4.4) and Proposition
2.4, we have\\

\noindent {\bf Proposition 4.3} {\it Let
$M=B\times_{b_1}F_1\times_{b_2}F_2\cdots
    \times_{b_m}F_m$ be a multiply twisted product and let $X,Y,Z\in \Gamma(TB)$
      and $V\in \Gamma(TF_i)$, $W\in \Gamma(TF_j)$, $U\in \Gamma(TF_k)$ and $P\in \Gamma(TB)$. Then}\\
\noindent $(1)\overline{R}(X,Y)Z=\overline{R}^B(X,Y)Z.$\\
\noindent $(2)\overline{R}(V,X)Y=-\left[\frac{H^{b_i}_B(X,Y)}{b_i}
+\frac{P(b_i)}{b_i}g(X,Y)+\pi(P)g(X,Y)+g(Y,\nabla_XP)-\pi(X)\pi(Y)\right]
V.$\\
\noindent $(3)\overline{R}(X,V)W=\overline{R}(V,W)X=\overline{R}(V,X)W=0 ~if ~i\neq j.$\\
\noindent $(4)\overline{R}(X,Y)V=0.$\\
\noindent $(5)\overline{R}(V,W)X=VX(lnb_i)W-WX(lnb_i)V ~ if ~i=j.$\\
\noindent $(6)\overline{R}(V,W)U=0~ if ~i=j\neq k~ or i\neq j \neq k.$\\
\noindent $(7)\overline{R}(U,V)W=-g(V,W)\frac{g_B({\rm
grad}_Bb_i,{\rm
grad}_Bb_k)}{b_ib_k}U-g(V,W)\left(\frac{P(b_i)}{b_i}+\frac{P(b_k)}{b_k}\right)U$\\
$~~~~~~~~~~-\pi(P)g(V,W)U, ~if ~i=j\neq k.$\\
\noindent
$(8)\overline{R}(X,V)W=[WX(lnb_i)]V-g(W,V)$\\
$\cdot\left[\frac{\nabla_X^B({\rm grad}_Bb_i)}{b_i}+\frac{{\rm
grad}_{F_i}(Xlnb_i)}{b_i^2}+\frac{P(b_i)}{b_i}X+\nabla_XP+\pi(P)X-\pi(X)P\right]
 ~if
~i=j.$\\
\noindent $(9)\overline{R}(U,V)W=g(U,W){\rm
grad}_B(V(lnb_i))-g(V,W){\rm
grad}_B(U(lnb_i))+R^{F_i}(U,V)W-\left(\frac{|{\rm
grad}_Bb_i|^2_B}{b_i^2}+2\frac{P(b_i)}{b_i}+\pi(P)\right)(g(V,W)U-g(U,W)V) ~if ~i=j=k.~$\\

\noindent {\bf Proposition 4.4} {\it Let
$M=B\times_{b_1}F_1\times_{b_2}F_2\cdots
    \times_{b_m}F_m$ be a multiply twisted product and let $X,Y,Z\in \Gamma(TB)$
      and $V\in \Gamma(TF_i)$, $W\in \Gamma(TF_j)$, $U\in \Gamma(TF_k)$ and $P\in \Gamma(TF_l)$. Then}\\
\noindent
$(1)\overline{R}(X,Y)Z={R}^B(X,Y)Z+\left[g(X,Z)\frac{Yb_l}{b_l}-g(Y,Z)\frac{Xb_l}{b_l}\right]P+\pi(P)[g(X,Z)Y-g(Y,Z)X]
.$\\
\noindent $(2)\overline{R}(V,X)Y=-\frac{H^{b_i}_B(X,Y)}{b_i}V
-\pi(P)g(X,Y)V~if ~i\neq l$\\
\noindent $(3)\overline{R}(V,X)Y=-\frac{H^{b_i}_B(X,Y)}{b_i}V
-\pi(V)\frac{Y(b_i)}{b_i}X-g(X,Y)\nabla_VP-g(X,Y)[\pi(P)V-\pi(V)P]
~if ~i= l$\\
\noindent $(4)\overline{R}(X,V)W=\frac{X(b_l)}{b_l}\pi(W)V~if ~i\neq j.$\\
\noindent
$(5)\overline{R}(V,W)X=-\delta^l_i\frac{\pi(V)}{b_i}X(b_i)W+\delta^l_j\frac{\pi(W)}{b_j}X(b_j)V~if ~i\neq j.$\\
\noindent $(6)\overline{R}(X,Y)V=\pi(V)\left[\frac{X(b_l)}{b_l}Y-\frac{Y(b_l)}{b_l}X\right].$\\
\noindent $(7)\overline{R}(V,W)X=VX(lnb_i)W-WX(lnb_i)V-\delta_i^l\frac{X(b_i)}{b_i}[\pi(V)W-\pi(W)V] ~ if ~i=j.$\\
\noindent $(8)\overline{R}(V,W)U=0~ if ~i=j\neq k~ or i\neq j \neq k.$\\
\noindent $(9)\overline{R}(U,V)W=-g(V,W)\frac{g_B({\rm
grad}_Bb_i,{\rm grad}_Bb_k)}{b_ib_k}U -g(W,\nabla_VP)U
-g(V,W)\nabla_UP-\pi(P)g(V,W)U
+g(V,W)\pi(U)P+\pi(W)[\pi(V)U-\pi(U)V]
, ~if ~i=j\neq k.$\\
\noindent
$(10)\overline{R}(X,V)W=[WX(lnb_i)]V-g(W,V)\frac{\nabla_X^B({\rm
grad}_Bb_i)}{b_i}-{\rm grad}_{F_i}(Xlnb_i)g_{F_i}(W,V)
+\frac{X(b_l)}{b_l}\pi(W)V-g(W,\nabla_VP)X-g(V,W)\frac{X(b_l)}{b_l}P-g(V,W)\pi(P)X+\pi(V)\pi(W)X
 ~if
~i=j.$\\
\noindent $(11)\overline{R}(U,V)W=g(U,W){\rm
grad}_B(V(lnb_i))-g(V,W){\rm
grad}_B(U(lnb_i))+R^{F_i}(U,V)W-\frac{|{\rm
grad}_Bb_i|^2_B}{b_i^2}(g(V,W)U-g(U,W)V)+\pi(P)[g(U,W)V-g(V,W)U]~if ~i=j=k\neq l.~$\\
\noindent $(12)\overline{R}(U,V)W=g(U,W){\rm
grad}_B(V(lnb_i))-g(V,W){\rm
grad}_B(U(lnb_i))+R^{F_i}(U,V)W-\frac{|{\rm
grad}_Bb_i|^2_B}{b_i^2}(g(V,W)U-g(U,W)V)
+g(W,\nabla_UP)V-g(W,\nabla_VP)U
+g(U,W)\nabla_VP-g(V,W)\nabla_UP+\pi(P)[g(U,W)V-g(V,W)U]
+[g(V,W)\pi(U)-g(U,W)\pi(V)]P+\pi(W)[\pi(V)U-\pi(U)V]
 ~if ~i=j=k=l.~$\\

By proposition 4.3 and 4.4, we have\\

\noindent {\bf Proposition 4.5} {\it Let
$M=B\times_{b_1}F_1\times_{b_2}F_2\cdots
    \times_{b_m}F_m$ be a multiply twisted product and let $X,Y,Z\in \Gamma(TB)$
      and $V\in \Gamma(TF_i)$, $W\in \Gamma(TF_j)$ and $P\in \Gamma(TB)$. Then}\\
\noindent $(1) \overline{{\rm Ric}} (X,Y)=\overline{{\rm
Ric}}^B(X,Y)+\sum_{i=1}^ml_i\left[ \frac{H_B^{b_i}(X,Y)}{b_i}
+\frac{P(b_i)}{b_i}g(X,Y)+\pi(P)g(X,Y)\right.$\\
$\left.~~~~~~~+g(Y,\nabla_XP)-\pi(X)\pi(Y)\right]
.$\\
\noindent $(2) \overline{{\rm Ric}} (X,V)=\overline{{\rm Ric}}
(V,X)=(l_i-1)[VX(lnb_i)].$\\
\noindent $(3) \overline{{\rm Ric}} (V,W)=0~if ~i\neq j.$\\
\noindent $(4) \overline{{\rm Ric}} (V,W)={\rm Ric}^{F_i}
(V,W)+\left[\frac{\triangle_Bb_i}{b_i}+(l_i-1)\frac{|{\rm
grad}_Bb_i|^2_B}{b_i^2}+\sum_{j\neq i}l_j\frac{g_B({\rm
grad}_Bb_i,{\rm grad}_Bb_j)}{b_ib_j}\right.$\\
$\left.+(\overline{n}-2)\pi(P)+\sum_{k=1}^n\varepsilon_k<\nabla_{E_k}P,E_k>+\sum_{j\neq
i}l_j\frac{Pb_j}{b_j}+(\overline{n}+l_i-2)\frac{Pb_i}{b_i}\right]
g(V,W)~if ~i= j,$\\
{\it where $E_k,~1\leq k\leq n$ is an orthonormal base of $B$ with
$\varepsilon_k=g(E_k,E_k)$ and ${\rm dim}B=n,~{\rm
dim}M=\overline{n}.$}\\

\noindent {\bf Corollary 4.6} {\it Let
$M=B\times_{b_1}F_1\times_{b_2}F_2\cdots
    \times_{b_m}F_m$ be a multiply twisted product and ${\rm
    dim}F_i>1$ and $P\in \Gamma(TB)$, then $(M,\overline{\nabla})$ is mixed Ricci-flat if and only if $M$ can
    be expressed as a multiply warped product. In particular, if $(M,\overline{\nabla})$
    is Einstein, then $M$ can be expressed as a multiply warped product.}\\

\noindent {\bf Proposition 4.7} {\it Let
$M=B\times_{b_1}F_1\times_{b_2}F_2\cdots
    \times_{b_m}F_m$ be a multiply twisted product and let $X,Y,Z\in \Gamma(TB)$
      and $V\in \Gamma(TF_i)$, $W\in \Gamma(TF_j)$ and $P\in \Gamma(TF_r)$. Then}\\
\noindent $(1) \overline{{\rm Ric}} (X,Y)={{\rm
Ric}}^B(X,Y)+\sum_{i=1}^ml_i\frac{H_B^{b_i}(X,Y)}{b_i}
+g(X,Y)\pi(P)(\overline{n}-2)$\\
$~~~~~+g(X,Y)\sum_{j_r=1}^{l_r}\varepsilon_{j_r}g(\nabla_{E_{j_r}^r}P,E_{j_r}^r)
.$\\
\noindent $(2) \overline{{\rm Ric}} (X,V)=(l_i-1)[VX(lnb_i)]+(\overline{n}-2)\frac{X(b_r)}{b_r}\pi(V).$\\
\noindent $(3) \overline{{\rm Ric}} (V,X)=(l_i-1)[VX(lnb_i)]+(2-\overline{n})\frac{X(b_r)}{b_r}\pi(V).$\\
\noindent $(3) \overline{{\rm Ric}} (V,W)=0~if ~i\neq j.$\\
\noindent $(4) \overline{{\rm Ric}} (V,W)={\rm Ric}^{F_i}
(V,W)+g(V,W)\left[\frac{\triangle_Bb_i}{b_i}+(l_i-1)\frac{|{\rm
grad}_Bb_i|^2_B}{b_i^2}+\sum_{j\neq i}l_j\frac{g_B({\rm
grad}_Bb_i,{\rm grad}_Bb_j)}{b_ib_j}\right.$\\
$\left.+(\overline{n}-2)\pi(P)\right]+(\overline{n}-2)g(W,\nabla_VP)+(2-\overline{n})\pi(V)\pi(W)+g(V,W){\rm
div}_{F_r}P~if ~i= j
,$\\

\noindent {\bf Corollary 4.8} {\it Let
$M=B\times_{b_1}F_1\times_{b_2}F_2\cdots
    \times_{b_m}F_m$ be a multiply twisted product and ${\rm
    dim}F_i>1$ and $P\in \Gamma(TF_r)$, then $(M,\overline{\nabla})$ is mixed Ricci-flat if and only if $M$ can
    be expressed as a multiply warped product and $b_r$ is only dependent on $F_r$. In particular, if $(M,\overline{\nabla})$
    is Einstein, then $M$ can be expressed as a multiply warped product.}\\

\noindent {\bf Proposition 4.9} {\it Let
$M=B\times_{b_1}F_1\times_{b_2}F_2\cdots
    \times_{b_m}F_m$ be a multiply twisted product and $P\in \Gamma(TB)$, then the scalar
    curvature $\overline{S}$ has the following expression:}\\
    $$
    \overline{S}=\overline{S}^B+2\sum_{i=1}^m\frac{l_i}{b_i}\triangle_Bb_i+\sum_{i=1}^m\frac{S^{F_i}}{b_i^2}+\sum_{i=1}^ml_i(l_i-1)
\frac{|{\rm grad}_Bb_i|^2_B}{b_i^2}$$ $$+\sum_{i=1}^m\sum_{j\neq
i}l_il_j\frac{g_B({\rm grad}_Bb_i,{\rm grad}_Bb_j)}{b_ib_j}+
\sum_{i=1}^ml_i(n+\overline{n}+l_i-2)\frac{P(b_i)}{b_i}$$
$$+\sum_{i=1}^m\sum_{j\neq
i}l_il_j\frac{P(b_j)}{b_j}+\sum_{i=1}^ml_i(n+\overline{n}-3)\pi(P)+2\sum_{i=1}^ml_i{\rm
div}_BP .\eqno(4.5)$$\\

\noindent {\bf Proposition 4.10} {\it Let
$M=B\times_{b_1}F_1\times_{b_2}F_2\cdots
    \times_{b_m}F_m$ be a multiply twisted product and $P\in \Gamma(TF_r)$, then the scalar
    curvature $\overline{S}$ has the following expression:}\\
    $$
    \overline{S}={S}^B+2\sum_{i=1}^m\frac{l_i}{b_i}\triangle_Bb_i+\sum_{i=1}^m\frac{S^{F_i}}{b_i^2}+\sum_{i=1}^ml_i(l_i-1)
\frac{|{\rm grad}_Bb_i|^2_B}{b_i^2}$$ $$+\sum_{i=1}^m\sum_{j\neq
i}l_il_j\frac{g_B({\rm grad}_Bb_i,{\rm grad}_Bb_j)}{b_ib_j}+
\pi(P)(\overline{n}-1)(\overline{n}-2)+2(\overline{n}-1){\rm div}_
{F_r}P.\eqno(4.6)$$\\

\noindent {\bf 4.2 Special multiply warped product with a
semi-symmetric connection}\\

 \indent Let $M=I\times_{b_1}F_1\times_{b_2}F_2\cdots
    \times_{b_m}F_m$ be a multiply warped product with the metric tensor $-dt^2\oplus  b_1^2g_{F_1}\oplus \cdots\oplus  b_m^2g_{F_m}$ and $I$ is an open
    interval in $\mathbb{R}$ and $b_i\in C^{\infty}(I)$.\\

\noindent {\bf Theorem 4.11} {\it  Let
$M=I\times_{b_1}F_1\times_{b_2}F_2\cdots
    \times_{b_m}F_m$ be a multiply warped product with the metric tensor $-dt^2\oplus b_1^2g_{F_1}\oplus \cdots\oplus
    b_m^2g_{F_m}$ and $P=\frac{\partial}{\partial t}$. Then
    $(M,\overline{\nabla})$ is Einstein with the Einstein constant
    $\lambda$ if and only if the following conditions are satisfied
    for any $i\in\{1,\cdots,m\}$}\\
\noindent (1){\it  $(F_i,\nabla^{F_i})$ is Einstein with the
Einstein
constant $\lambda_i$, $i\in\{1,\cdots,m\}$.}\\
\noindent (2)
$\sum_{i=1}^ml_i\left(\frac{b_i'}{b_i}-\frac{b_i''}{b_i}\right)=\lambda.$\\
\noindent
(3)$\lambda_i-b_ib_i''-(l_i-1)b_i'^2+(b_i^2-b_ib_i')\sum_{j\neq
i}l_j\frac{b_j'}{b_j}+(2-\overline{n})b_i^2+(\overline{n}+l_i-2)b_ib_i'=\lambda
b_i^2.$\\

\noindent{\bf Proof.} By Proposition 4.5, we have
$$\overline{{\rm Ric}}\left(\frac{\partial}{\partial t},\frac{\partial}{\partial
t}\right)=-\sum_{i=1}^ml_i\left(\frac{b_i'}{b_i}-\frac{b_i''}{b_i}\right);\eqno(4.7)$$
$$\overline{{\rm Ric}}\left(\frac{\partial}{\partial t},V\right)=\overline{{\rm Ric}}\left(V,\frac{\partial}{\partial
t}\right)=0;\eqno(4.8)$$
$$\overline{{\rm Ric}}\left(V,W\right)={\rm
Ric}^{F_i}(V,W)+g_{F_i}(V,W)\left[-b_ib_i"-(l_i-1)b_i'^2+(b_i^2-b_ib_i')\sum_{j\neq
i}l_j\frac{b_j'}{b_j}\right.$$
$$\left.+(2-\overline{n})b_i^2+(\overline{n}+l_i-2)b_ib_i'\right].\eqno(4.9)$$
By (4.7)-(4.9) and the Einstein condition, we get the above
theorem.~~~~ $\Box$\\

\noindent {\bf Theorem 4.12} {\it  Let
$M=I\times_{b_1}F_1\times_{b_2}F_2\cdots
    \times_{b_m}F_m$ be a multiply warped product with the metric tensor $-dt^2\oplus b_1^2g_{F_1}\oplus \cdots\oplus
    b_m^2g_{F_m}$ and $P\in \Gamma(TF_r)$ with $g_{F_r}(P,P)=1$ and $\overline{n}>2$. Then
    $(M,\overline{\nabla})$ is Einstein with the Einstein constant
    $\lambda$ if and only if the following conditions are satisfied
    for any $i\in\{1,\cdots,m\}$}\\
\noindent (1){\it  $(F_i,\nabla^{F_i})~(i\neq r)$ is Einstein with
the Einstein
constant $\lambda_i$, $i\in\{1,\cdots,m\}$.}\\
\noindent (2){\it $b_r$ is a constant and
$\sum_{i=1}^ml_i\frac{b_i''}{b_i}=\mu_0;~{\rm div}_{F_r}P=\mu_1,~\mu_0-\mu_1+\lambda=(2-\overline{n})b_r^2,$ where $\mu_0,\mu_1$ are constants.}\\
\noindent (3){\it $ {\rm
Ric}^{F_r}(V,W)+\overline{\lambda}g_{F_r}(V,W)=(\overline{n}-2)\left[\pi(V)\pi(W)-g(W,\nabla_VP)\right],~for~
V,W\in \Gamma(TF_r).$}\\
 \noindent
(4){\it
$\lambda_i-b_ib_i''+(\overline{n}-2)b_i^2b_r^2-b_ib_i'\sum_{j\neq
i}l_j\frac{b_j'}{b_j}-(l_i-1)b_i'^2=(\lambda-\mu_1)b_i^2.$}\\

\noindent{\bf Proof.} By Proposition 4.7 (2) and $g_{F_r}(P,P)=1$,
we have $b_r$ is a constant. By Proposition 4.7, then
,$$\overline{{\rm Ric}}\left(\frac{\partial}{\partial
t},\frac{\partial}{\partial
t}\right)=\sum_{i=1}^ml_i\frac{b_i''}{b_i}+(2-\overline{n})b_r^2-{\rm
div}_{F_r}P=-\lambda;\eqno(4.10)$$ By variables separation, we have
$$\sum_{i=1}^ml_i\frac{b_i''}{b_i}=\mu_0,~{\rm div}_{F_r}P=\mu_1.~\mu_0-\mu_1+\lambda=(2-\overline{n})b_r^2,
.\eqno(4.11)$$

$$  \overline{{\rm Ric}} (V,W)={\rm Ric}^{F_i}
(V,W)+b_i^2g_{F_i}(V,W)\left[-\frac{b_i''}{b_i}+(l_i-1)\frac{-b_i'^2}{b_i^2}+\sum_{j\neq
i}l_j\frac{-b_i'b_j'}{b_ib_j}\right.$$
$$\left.+(\overline{n}-2)\pi(P)\right]+(\overline{n}-2)g(W,\nabla_VP)+(2-\overline{n})\pi(V)\pi(W)+g(V,W){\rm
div}_{F_r}P.\eqno(4.12)$$ When $i\neq r$, then $\nabla_VP=\pi(V)=0$,
so
$$  \overline{{\rm Ric}} (V,W)={\rm Ric}^{F_i}
(V,W)+b_i^2g_{F_i}(V,W)\left[-\frac{b_i''}{b_i}+(l_i-1)\frac{-b_i'^2}{b_i^2}+\sum_{j\neq
i}l_j\frac{-b_i'b_j'}{b_ib_j}\right.$$
$$\left.+(\overline{n}-2)b_r^2\right]+\mu_1b_i^2g_{F_i}(V,W)=\lambda
b_i^2g_{F_i}(V,W).\eqno(4.13)$$ By variables separation, we have
$(F_i,\nabla^{F_i})~(i\neq r)$ is Einstein with the Einstein
constant $\lambda_i$ and
$$\lambda_i-b_ib_i''+(\overline{n}-2)b_i^2b_r^2-b_ib_i'\sum_{j\neq
i}l_j\frac{b_j'}{b_j}-(l_i-1)b_i'^2=(\lambda-\mu_1)b_i^2.\eqno(4.14)$$
When $i=r$ and $b_r$ is a constant, then
$${\rm Ric}^{F_i}
(V,W)+b_r^2[(\overline{n}-2)b_r^2+\mu_1-\lambda]g_{F_i}(V,W)=(\overline{n}-2)\left[\pi(V)\pi(W)-g(W,\nabla_VP)\right].\eqno(4.15)$$
So we prove the above theorem. ~~~~~$\Box$\\

When $M=I\times_{b_1}F_1\times_{b_2}F_2\cdots
    \times_{b_m}F_m$ be a multiply warped product and $P=\frac{\partial}{\partial
    t}$, by Proposition 4.9, we have
 $$
    \overline{S}=-2\sum_{i=1}^ml_i\frac{b_i''}{b_i}+\sum_{i=1}^m\frac{S^{F_i}}{b_i^2}+\sum_{i=1}^ml_i(l_i-1)
\frac{-b_i'^2}{b_i^2}+\sum_{i=1}^m\sum_{j\neq
i}l_il_j\frac{-b_ib_j}{b_ib_j}$$ $$+
\sum_{i=1}^ml_i(\overline{n}+l_i-1)\frac{b_i'}{b_i}+\sum_{i=1}^m\sum_{j\neq
i}l_il_j\frac{b_j'}{b_j}-\sum_{i=1}^ml_i(\overline{n}-2) .\eqno(4.16)$$\\
The following result just follows from the method of separation of
variables and the fact that each $S^{F_i}$ is function defined on
$F_i$.\\

\noindent {\bf Proposition 4.13} {\it  Let
$M=I\times_{b_1}F_1\times_{b_2}F_2\cdots
    \times_{b_m}F_m$ be a multiply warped product and $P=\frac{\partial}{\partial
    t}$. If
    $(M,\overline{\nabla})$ has constant scalar curvature
    $\overline{S}$, then each $(F_i,\nabla^{F_i})$ has constant
    scalar curvature $S^{F_i}$.}\\

When $P\in\Gamma(TF_r)$, by Proposition 4.10, we have
$$
    \overline{S}=-2\sum_{i=1}^ml_i\frac{b_i''}{b_i}+\sum_{i=1}^m\frac{S^{F_i}}{b_i^2}+\sum_{i=1}^ml_i(l_i-1)
\frac{-b_i'^2}{b_i^2}$$ $$+\sum_{i=1}^m\sum_{j\neq
i}l_il_j\frac{-b_i'b_j'}{b_ib_j}+
\pi(P)(\overline{n}-1)(\overline{n}-2)+2(\overline{n}-1){\rm div}_
{F_r}P.\eqno(4.17)$$\\

\noindent {\bf Proposition 4.14} {\it  Let
$M=I\times_{b_1}F_1\times_{b_2}F_2\cdots
    \times_{b_m}F_m$ be a multiply warped product and $P\in\Gamma(TF_r)$. If
    $(M,\overline{\nabla})$ has constant scalar curvature
    $\overline{S}$, then each $(F_i,\nabla^{F_i})~(i\neq r)$ has constant
    scalar curvature $S^{F_i}$ and if $g_{F_r}(P,P)$ and ${\rm div}_
{F_r}P$ are constants, then $S^{F_r}$ is also a constant. }\\

\noindent{\bf 4.3 Generalized Robertson-Walker spacetimes with a
semi-symmetric metric connection}\\

\indent In this section, we study $M=I\times F$ with the metric
tensor $-dt^2+f(t)^2g_F$. As a corollary of Theorem 4.11, we
obtain:\\

\noindent {\bf Corollary 4.15} {\it Let $M=I\times F$ with the
metric tensor $-dt^2+f(t)^2g_F$ and $P=\frac{\partial}{\partial t}$.
Then
    $(M,\overline{\nabla})$ is Einstein with the Einstein constant
    $\lambda$ if and only if the following conditions are satisfied
    }\\
\noindent (1){\it  $(F,\nabla^{F})$ is Einstein with the Einstein
constant $\lambda_F$.}\\
\noindent (2)
$l\left(\frac{f'}{f}-\frac{f''}{f}\right)=\lambda.$\\
\noindent
(3)$\lambda_F-ff"+(1-l)f'^2+(1-l-\lambda)f^2+(2l-1)f'f=0.$\\

\noindent {\bf Remark.} In Theorem 5.1 in [SO], they got the
Einstein condition of $M=I\times F$ with a semi-symmetric metric
connection, but they did not consider the above conditions (2) and
(3).\\

\noindent {\bf Corollary 4.16} {\it Let $M=I\times F$ with the
metric tensor $-dt^2+f(t)^2g_F$ and $P=\frac{\partial}{\partial t}$
and ${\rm dim}F=1$. Then  $(M,\overline{\nabla})$ is Einstein with
the Einstein constant
    $\lambda$ if and only if $f''=f'-\lambda f$.}\\

By Corollary 4.15 (2) and (3), we get\\

\noindent {\bf Corollary 4.17} {\it Let $M=I\times F$ with the
metric tensor $-dt^2+f(t)^2g_F$ and $P=\frac{\partial}{\partial t}$
and ${\rm dim}F>1$. Then
    $(M,\overline{\nabla})$ is Einstein with the Einstein constant
    $\lambda$ if and only if the following conditions are satisfied
    }\\
\noindent (1){\it  $(F,\nabla^{F})$ is Einstein with the Einstein
constant $\lambda_F$.}\\
\noindent (2)
$f''=f'-\frac{\lambda}{l}f.$\\
\noindent
(3)$\frac{\lambda_F}{1-l}+f'^2+(1+\frac{\lambda}{l})f^2-2ff'=0.$\\

By Corollary 4.16 and elementary methods for ordinary differential
 equations, we get\\

 \noindent {\bf Theorem 4.18} {\it Let $M=I\times F$ with the
metric tensor $-dt^2+f(t)^2g_F$ and $P=\frac{\partial}{\partial t}$
and ${\rm dim}F=1$. Then  $(M,\overline{\nabla})$ is Einstein with
the Einstein constant
    $\lambda$ if and only if }\\
\noindent (1)
$\lambda<\frac{1}{4},~~f(t)=c_1e^{\frac{1+\sqrt{1-4\lambda}}{2}t}+c_2e^{\frac{1-\sqrt{1-4\lambda}}{2}t},$\\
\noindent (2)
$\lambda=\frac{1}{4},~~f(t)=c_1e^{\frac{1}{2}t}+c_2te^{\frac{1}{2}t},$\\
\noindent (3) $\lambda>\frac{1}{4},~~f(t)=c_1e^{\frac{1}{2}t}{\rm
cos}\left(\frac{\sqrt{4\lambda-1}}{2}t\right)+
c_2e^{\frac{1}{2}t}{\rm sin}\left(\frac{\sqrt{4\lambda-1}}{2}t\right),$\\

Let $\frac{\lambda}{l}=d_0$,~$\frac{\lambda_F}{1-l}=\overline{d_0}$,
$\frac {1+\sqrt{1-4d_0}}{2}=a_0,~~\frac {1-\sqrt{1-4d_0}}{2}=b_0, $ then $a_0+b_0=1, d_0=a_0b_0$. When ${\rm dim}F>1$, by Corollary 4.17 (2)\\

\noindent {\bf Case i)} $d_0<\frac{1}{4}$, then
$f=c_1e^{a_0t}+c_2e^{b_0t}.$ By Corollary 4.17 (3), then
$$\overline{d_0}+c_1^2(a_0^2+1+a_0b_0-2a_0)e^{2a_0t}+c_2^2(b_0^2+1+a_0b_0-2b_0)e^{2b_0t}$$
$$+2c_1c_2(2a_0b_0+1-a_0-b_0)e^{(a_0+b_0)t}=0 .\eqno(4.18)$$
When $b_0=0$, we get $d_0=0,~ a_0=1,\lambda=0$. By (4.18),
$\overline{d_0}+c_2^2=0$, so $\lambda_F=(l-1)c_2^2.$ In this case
$f=c_1e^t+c_2.$ When $b\neq 0$, then $e^{2a_0t}$, $e^{2b_0t}$ and
$e^{(a_0+b_0)t}$ are linear independent, so
$c_2^2(b_0^2+1+a_0b_0-2b_0)=c_2^2(1-b_0)=0$ and $c_2=0$. Then
$c_1\neq 0$, by $c_1^2(a_0^2+1+a_0b_0-2a_0)=c_1^2(1-a_0)=0$, so
$a_0=1$, then $d_0=\lambda=0$. Thus $f=c_1e^t$.\\

\noindent {\bf Case ii)} $d_0=\frac{1}{4}$, then
$f=c_1e^{\frac{1}{2}t}+c_2te^{\frac{1}{2}t}.$ By Corollary 4.17 (3),
then
$$\overline{d_0}
+e^t\left[(\frac{1}{2}c_1+c_2+\frac{1}{2}c_2t)^2+\frac{5}{4}(c_1+c_2t)^2-2(\frac{1}{2}c_1+c_2+\frac{1}{2}c_2t)(c_1+c_2t)\right]=0.\eqno(4.19)$$
The coefficient of $t^2e^t$ is $\frac{1}{2}c_2^2$, so $c_2=0$. The
coefficient of $e^t$ is
$(\frac{1}{2}c_1+c_2)^2+\frac{5}{4}(c_1)^2-2c_1(\frac{1}{2}c_1+c_2)$,
so $c_1=0$, in this case we have no solutions.\\

\noindent {\bf Case iii)} $d_0>\frac{1}{4}$, then
$f(t)=c_1e^{\frac{1}{2}t}{\rm cos}(h_0t)+c_2e^{\frac{1}{2}t}{\rm
sin}(h_0t),$ where $h_0=\frac{\sqrt{4d_0-1}}{2}.$  By Corollary 4.17
(3), then
$$\overline{d_0}
+e^t\left\{\left[(\frac{c_1}{2}+c_2h_0){\rm
cos}(h_0t)+(\frac{c_2}{2}-c_1h_0){\rm sin}(h_0t)\right]^2\right.$$
$$+(1+d_0)(c_1{\rm cos}(h_0t)+c_2{\rm sin}(h_0t))^2$$
$$\left.-2
(c_1{\rm cos}(h_0t)+c_2{\rm
sin}(h_0t))\left[(\frac{c_1}{2}+c_2h_0){\rm
cos}(h_0t)+(\frac{c_2}{2}-c_1h_0){\rm
sin}(h_0t)\right]\right\}=0.\eqno(4.20)$$ Consider the coefficients
of ${\rm cos}^2(h_0t)e^t$ and ${\rm sin}^2(h_0t)e^t$, we get
$$(\frac{1}{4}+d_0)c_1^2+c_2^2h_0^2-c_1c_2h_0=0;~~(\frac{1}{4}+d_0)c_2^2+c_1^2h_0^2+c_1c_2h_0=0.\eqno(4.21)$$
Plusing the above two equalities ,then $\frac{1}{4}+d_0+h_0^2=0$ and
$d_0=0$. There is a contradiction with $d_0>\frac{1}{4}$ and in this
case we have no solutions. So we obtain the following theorem.\\

\noindent {\bf Theorem 4.19} {\it Let $M=I\times F$ with the metric
tensor $-dt^2+f(t)^2g_F$ and $P=\frac{\partial}{\partial t}$ and
${\rm dim}F>1$. Then
    $(M,\overline{\nabla})$ is Einstein with the Einstein constant
    $\lambda$ if and only if $\lambda=0$ and $f=c_1e^t+c_2$ and
 $(F,\nabla^{F})$ is Einstein with the Einstein
constant $(l-1)c_2^2$.}\\

\indent By (4.16) and (4.17), we have\\

\noindent {\bf Corollary 4.20} {\it  Let  Let $M=I\times F$ with the
metric tensor $-dt^2+f(t)^2g_F$ and $P=\frac{\partial}{\partial t}$.
 If
    $(M,\overline{\nabla})$ has constant scalar curvature
    $\overline{S}$ if and only if $(F,\nabla^{F})$ has constant
    scalar curvature $S^{F}$ and}
$$\overline{S}=\frac{S^F}{f^2}-2l\frac{f''}{f}-l(l-1)\frac{f'^2}{f^2}+2l^2\frac{f'}{f}+(1-l)l.\eqno(4.22)$$\\

\noindent {\bf Corollary 4.21} {\it  Let  Let $M=I\times F$ with the
metric tensor $-dt^2+f(t)^2g_F$ and $P\in\Gamma(TF)$ and
$g_F(P,P)=c_0, ~{\rm div}_{F_r}P=c_0'$.
 If
    $(M,\overline{\nabla})$ has constant scalar curvature
    $\overline{S}$ if and only if $(F,\nabla^{F})$ has constant
    scalar curvature $S^{F}$ and}
$$\overline{S}=\frac{S^F}{f^2}-2l\frac{f''}{f}-l(l-1)\frac{f'^2}{f^2}+c_0(l-1)lf^2+2c_0'l.\eqno(4.23)$$

\indent In (4.22), we make the change of variable $f(t)=\sqrt{v(t)}$
and have the following equation
$$v''(t)+\frac{l-3}{4}\frac{v'(t)^2}{v(t)}-lv'(t)+(l-1+\frac{\overline{S}}{l})v(t)-\frac{S^F}{l}=0.\eqno(4.24)$$\\

\noindent {\bf Theorem 4.22} {\it  Let $M=I\times F$ with the metric
tensor $-dt^2+f(t)^2g_F$ and $P=\frac{\partial}{\partial t}$ and
${\rm dim} F=l=3$.
 If
    $(M,\overline{\nabla})$ has constant scalar curvature
    $\overline{S}$ if and only if $(F,\nabla^{F})$ has constant
    scalar curvature $S^{F}$ and}\\
    \noindent (1) $\overline{S}<\frac{3}{4}$ and $\overline{S}\neq
    -6$,~~$v(t)=c_1e^{\frac{3+\sqrt{1-\frac{4}{3}\overline{S}}}{2}t}+c_2e^{\frac{3-\sqrt{1-\frac{4}{3}\overline{S}}}{2}t}+\frac{S^F}{6+\overline{S}}.$\\
\noindent (2) $\overline{S}=\frac{3}{4},~~
~~v(t)=c_1e^{\frac{3}{2}t}+c_2te^{\frac{3}{2}t}+\frac{S^F}{6+\overline{S}}.$\\
\noindent (3) $\overline{S}>\frac{3}{4},~~
~~v(t)=c_1e^{\frac{3}{2}t}{\rm
cos}\left(\frac{\sqrt{\frac{4}{3}\overline{S}-1}}{2}t\right)
+c_2e^{\frac{3}{2}t}{\rm
sin}\left(\frac{\sqrt{\frac{4}{3}\overline{S}-1}}{2}t\right)+\frac{S^F}{6+\overline{S}}.$\\
\noindent (4)
$\overline{S}=-6,~~~~v(t)=c_1-\frac{S^F}{9}t+c_2e^{3t}.$\\

\noindent{\bf Proof.} If $l=3$, then we have a simple differential
equation
$$v''(t)-3v'(t)+(2+\frac{\overline{S}}{3})v(t)-\frac{S^F}{3}=0.\eqno(4.25)$$\\
If $\overline{S}\neq
    -6$, we putting
    $h(t)=(2+\frac{\overline{S}}{3})v(t)-\frac{S^F}{3},$ it follows
    that $h''(t)-3h'(t)+(2+\frac{\overline{S}}{3})h(t)=0$. The above
    solutions (1)-(3) follow directly from elementary methods for
    ordinary differential equations. When $\overline{S}=
    -6$, then $v''(t)-3v'(t)-\frac{S^F}{3}=0$, we get the solution
    (4).~~~~~~~~~~~~~$\Box$\\

\noindent {\bf Theorem 4.23} {\it  Let $M=I\times F$ with the metric
tensor $-dt^2+f(t)^2g_F$ and $P=\frac{\partial}{\partial t}$ and
${\rm dim} F=l\neq 3$ and $S^F=0$.
 If
    $(M,\overline{\nabla})$ has constant scalar curvature
    $\overline{S}$ if and only if }\\
    \noindent (1)
    $\overline{S}<\frac{l}{l+1},~~v(t)=\left(c_1e^{\frac{l+\sqrt{1-\frac{l+1}{l}\overline{S}}}{2}t}
    +c_2e^{\frac{l-\sqrt{1-\frac{l+1}{l}\overline{S}}}{2}t}\right)^{\frac{4}{l+1}}.$\\

 \noindent (2)
    $\overline{S}=\frac{l}{l+1},~~v(t)=\left(c_1e^{\frac{l}{2}t}
    +c_2te^{\frac{l}{2}t}\right)^{\frac{4}{l+1}}.$\\

\noindent (3)
    $\overline{S}>\frac{l}{l+1},~~v(t)=\left(c_1e^{\frac{l}{2}t}{\rm
    cos}\left(\frac{\sqrt{\frac{l+1}{l}\overline{S}-1}}{2}t\right)
    +c_2e^{\frac{l}{2}t}{\rm
    sin}\left(\frac{\sqrt{\frac{l+1}{l}\overline{S}-1}}{2}t\right)\right)^{\frac{4}{l+1}}.$\\

\noindent{\bf Proof.} In this case, the equation (4.24) is changed
into the simpler form
$$\frac{v''(t)}{v(t)}+\frac{l-3}{4}\frac{v'(t)^2}{v(t)^2}-l\frac{v'(t)}{v(t)}+(l-1+\frac{\overline{S}}{l})=0.\eqno(4.26)$$\\
Putting $v(t)=w(t)^{\frac{4}{l+1}}$, then $w(t)$ satisfies the
equation $w''-lw'+\frac{(l+1)}{4}(l-1+\frac{\overline{S}}{l})w=0$,
by the elementary methods for
    ordinary differential equations, we prove the above
    theorem.~~~~~~~~~$\Box$\\

When ${\rm dim} F=l\neq 3$ and $S^F\neq 0$, putting
$v(t)=w(t)^{\frac{4}{l+1}}$, then $w(t)$ satisfies the equation
$$w''-lw'+\frac{(l+1)}{4}(l-1+\frac{\overline{S}}{l})w-\frac{(l+1)}{4}\frac{S^F}{l}w^{1-\frac{4}{l+1}}=0.\eqno(4.27)$$\\

\noindent {\bf 4.4 ~Generalized Kasner spacetimes with a
semi-symmetric metric
 connection }\\

 \indent In this section, we consider the scalar and Ricci
 curvature of generalized Kasner spacetimes with a semi-symmetric metric
 connection. We recall the definition of generalized Kasner
 spacetimes ([DU1]).\\

 \noindent{\bf Definition 4.24} A generalized Kasner spacetime
 $(M,g)$ is a Lorentzian multiply warped product of the form
 $M=I\times_{\phi^{p_1}}F_1\times\cdots\times_{\phi^{p_m}}F_m$ with
 the metric
 $g=-dt^2\oplus\phi^{2p_1}g_{F_1}\oplus\cdots\oplus\phi^{2p_m}g_{F_m}$,
 where $\phi:~I\rightarrow (0,\infty)$ is smooth and
 $p_i\in\mathbb{R}$, for any $i\in\{1,\cdots,m\}$ and also
 $I=(t_1,t_2)$.\\

\indent We introduce the following parameters
$\zeta=\sum_{i=1}^ml_ip_i$ and $\eta=\sum_{i=1}^ml_ip_i^2$ for
generalized Kasner spacetimes. By Theorem 4.11 and direct
computations, we get\\

\noindent {\bf Proposition 4.25}{\it~Let
$M=I\times_{\phi^{p_1}}F_1\times\cdots\times_{\phi^{p_m}}F_m$ be a
generalized Kasner spacetime and $P=\frac{\partial}{\partial t}$.
Then
    $(M,\overline{\nabla})$ is Einstein with the Einstein constant
    $\lambda$ if and only if the following conditions are satisfied
    for any $i\in\{1,\cdots,m\}$}\\
\noindent (1){\it  $(F_i,\nabla^{F_i})$ is Einstein with the
Einstein
constant $\lambda_i$, $i\in\{1,\cdots,m\}$.}\\
\noindent (2)
$\zeta\left(\frac{\phi'-\phi''}{\phi}\right)-(\eta-\zeta)\frac{\phi'^2}{\phi^2}=\lambda.$\\
\noindent
(3)$\frac{\lambda_i}{\phi^{2p_i}}-p_i\frac{\phi''}{\phi}-(\zeta-1)p_i\frac{\phi'^2}{\phi^2}+[\zeta+(\overline{n}-2)p_i]\frac{\phi'}{\phi}=
\overline{n}+\lambda-2.$\\

By (4.16) we obtain\\

\noindent {\bf Proposition 4.26}{\it~Let
$M=I\times_{\phi^{p_1}}F_1\times\cdots\times_{\phi^{p_m}}F_m$ be a
generalized Kasner spacetime and $P=\frac{\partial}{\partial t}$.
Then
    $(M,\overline{\nabla})$ has constant scalar curvature
    $\overline{S}$ if and only if each $(F_i,\nabla^{F_i})$ has constant
    scalar curvature $S^{F_i}$ and}
$$\overline{S}=\sum_{i=1}^m\frac{S^{F_i}}{\phi^{2p_i}}-2\zeta\frac{\phi''}{\phi}-(\eta+\zeta^2-2\zeta)\frac{\phi'^2}{\phi^2}+2(\overline{n}-1)\zeta
\frac{\phi'}{\phi}+(2-\overline{n})(\overline{n}-1).\eqno(4.28)$$

\indent Nextly, we first give a classification of four-dimensional
generalized Kasner spacetimes with a semi-symmetric metric
 connection and then consider Ricci tensors and scalar curvatures
 of them.\\

\noindent{\bf Definition 4.27} Let
 $M=I\times_{b_1}F_1\times\cdots\times_{b_m}F_m$ with
 the metric
 $g=-dt^2\oplus b_1^2g_{F_1}\oplus\cdots\oplus b_m^2g_{F_m}$.\\
\noindent {\bf $\cdot$} $(M,g)$ is said to be of Type (I) if $m=1$
and ${\rm dim}(F)=3$.\\
\noindent {\bf $\cdot$} $(M,g)$ is said to be of Type (II) if $m=2$
and ${\rm dim}(F_1)=1$ and ${\rm dim}(F_2)=2$.\\
\noindent {\bf $\cdot$} $(M,g)$ is said to be of Type (III) if $m=3$
and ${\rm dim}(F_1)=1$,~${\rm dim}(F_2)=1$ and ${\rm dim}(F_3)=1$.\\

By Theorem 4.19 and 4.22, we have given a classification of Type (I)
Einstein spaces and Type (I) spaces with the constant scalar
curvature.\\

\noindent{\bf $\cdot$ Classification of Einstein Type (II)
generalized Kasner space-times with a semi-symmetric metric
 connection}\\
 \indent Let $M=I\times_{\phi^{p_1}}F_1\times_{\phi^{p_2}}F_2$ be an
 Einstein type (II) generalized Kasner spacetime and $P=\frac{\partial}{\partial
 t}$. Then $\zeta=p_1+2p_2$, $\eta=p_1^2+2p_2^2$.  By Proposition
 (4.25), we have
          $$\zeta\left(\frac{\phi'-\phi''}{\phi}\right)-(\eta-\zeta)\frac{\phi'^2}{\phi^2}=\lambda,\eqno(4.29i)$$
        $$ -p_1\frac{\phi''}{\phi}-(\zeta-1)p_1\frac{\phi'^2}{\phi^2}+[\zeta+2p_1]\frac{\phi'}{\phi}=
\lambda+2,\eqno(4.29ii)$$
       $$ \frac{\lambda_2}{\phi^{2p_2}}-p_2\frac{\phi''}{\phi}-(\zeta-1)p_2\frac{\phi'^2}{\phi^2}
        +[\zeta+2p_2]\frac{\phi'}{\phi}=
\lambda+2,\eqno(4.29iii)$$ where $\lambda_2$ is a constant. Consider
following two cases:\\

\noindent {\bf Case i)}~ $\underline{\zeta=0}$\\
 \indent In this
case, $p_2=-\frac{1}{2}p_1$, $\eta=\frac{3}{2}p_1^2$. Then by
(4.29), we have
$$-\eta\frac{\phi'^2}{\phi^2}=\lambda ,\eqno(4.30i)$$
        $$ p_1\left(-\frac{\phi''}{\phi}+\frac{\phi'^2}{\phi^2}+2\frac{\phi'}{\phi}\right)=
\lambda+2,\eqno(4.30ii)$$
       $$ \frac{\lambda_2}{\phi^{-p_1}}-\frac{1}{2}p_1\left(-\frac{\phi''}{\phi}+\frac{\phi'^2}{\phi^2}+2\frac{\phi'}{\phi}\right)
      =
\lambda+2,\eqno(4.30iii)$$
 {\bf Case i a)}~ $\underline{\eta=0}$,\\
then $p_i=0$, by (4.30i), $\lambda=0$. By (4.30ii), $\lambda+2=0$,
this is a contradiction.\\
{\bf Case i b)}~ $\underline{\eta\neq 0}$,\\
then $p_i\neq 0$.\\
 {\bf Case i b)1)} $\underline{\lambda_2=0}$\\
by (4.30ii) and (4.30iii), $\lambda=-2$ and
$$-\frac{\phi''}{\phi}+\frac{\phi'^2}{\phi^2}+2\frac{\phi'}{\phi}=0,~~
\frac{\phi'^2}{\phi^2}=\frac{2}{\eta},\eqno(4.31)$$ then
$\phi=c_0e^{\pm\sqrt{\frac{2}{\eta}}t}$ which does not satisfy the
first equation in (4.31), this a contradiction.\\
 {\bf Case i b)2)} $\underline{\lambda_2\neq 0}$\\
by (4.30ii) and (4.30iii), we have
$\frac{\lambda_2}{\phi^{-p_1}}=\frac{3}{2}(\lambda+2),$ so $\phi$ is
a constant. By (4.30ii), $\lambda+2=0$, so $\lambda_2=0$, this is a
contradiction. In a word, we have no solutions when $\zeta=0.$\\

\noindent {\bf Case ii)}~ $\underline{\zeta\neq  0},$\\
then $\eta\neq 0$. Putting $\phi=\psi^{\frac{\zeta}{\eta}}$, then
$\psi''-\psi'+\frac{\lambda\eta}{\zeta^2}\psi=0$. Hence,\\
\noindent (1)~~$\lambda<\frac{\zeta^2}{4\eta},$
$\psi=c_1e^{\frac{1+\sqrt{1-\frac{4\lambda\eta}{\zeta^2}}}{2}t} +
c_2e^{\frac{1-\sqrt{1-\frac{4\lambda\eta}{\zeta^2}}}{2}t},$\\
\noindent (2)~~$\lambda= \frac{\zeta^2}{4\eta},$
$\psi=c_1e^{\frac{1}{2}t} + c_2te^{\frac{1}{2}t},$\\
\noindent (3)~~$\lambda>\frac{\zeta^2}{4\eta},$
$\psi=c_1e^{\frac{1}{2}t}{\rm cos}\left(
\frac{\sqrt{\frac{4\lambda\eta}{\zeta^2}-1}}{2}t\right)
 + c_2e^{\frac{1}{2}t}{\rm sin}\left(
\frac{\sqrt{\frac{4\lambda\eta}{\zeta^2}-1}}{2}t\right).$\\
We make (4.30) into
$$\frac{\zeta^2}{\eta}\frac{\psi'-\psi''}{\psi}=\lambda,~~\psi=\phi^{\frac{\eta}{\zeta}},\eqno(4.32i)$$
$$-\frac{p_1}{\zeta}\frac{(\phi^\zeta)''}{\phi^\zeta}+\frac{\zeta+2p_1}{\zeta}
\frac{(\phi^\zeta)'}{\phi^\zeta}=\lambda+2,\eqno(4.32ii)$$
$$\frac{\lambda_2}{\phi^{2p_2}}-\frac{p_2}{\zeta}\frac{(\phi^\zeta)''}{\phi^\zeta}+\frac{\zeta+2p_2}{\zeta}
\frac{(\phi^\zeta)'}{\phi^\zeta}=\lambda+2,\eqno(4.32iii)$$ When
$p_1=p_2$, the type (II) spaces turns into type (I) spaces, so we
assume $p_1\neq p_2$. By (4.32ii) and (4.32iii), then
$$\psi'=\frac{p_1\lambda_2\eta}{(p_2-p_1)\zeta^2}\psi^{1-\frac{2p_2\zeta}{\eta}}
+ \frac{(\lambda+2)\eta}{\zeta^2}\psi.\eqno(4.33)$$

\noindent {\bf Case ii)(1)}~~$\lambda<\frac{\zeta^2}{4\eta},$
$\psi=c_1e^{at} + c_2e^{bt},$ where
$a=\frac{1+\sqrt{1-\frac{4\lambda\eta}{\zeta^2}}}{2},
,~b=\frac{1-\sqrt{1-\frac{4\lambda\eta}{\zeta^2}}}{2}t.$\\
By (4.33),
$$ac_1e^{at} + bc_2e^{bt}=\frac{p_1\lambda_2\eta}{(p_2-p_1)\zeta^2}(c_1e^{at} + c_2e^{bt})^{1-\frac{2p_2\zeta}{\eta}}
+ \frac{(\lambda+2)\eta}{\zeta^2}(c_1e^{at} +
c_2e^{bt}).\eqno(4.34)$$
 \noindent {\bf Case ii)(1)(a)}~$\underline{c_1=0}$,\\
 then $$\left[b-\frac{(\lambda+2)\eta}{\zeta^2}\right]c_2e^{bt}=
\frac{p_1\lambda_2\eta}{(p_2-p_1)\zeta^2}(
c_2e^{bt})^{1-\frac{2p_2\zeta}{\eta}}.\eqno(4.35)$$
 \noindent {\bf Case ii)(1)(a)1)}~$\underline{b\neq 0},~\underline{p_1\lambda_2\neq
 0}$\\
 then $p_2=0$ and $\zeta=p_1,~\eta=p_1^2$ and
 $b=\frac{1-\sqrt{1-4\lambda}}{2}$ and $\psi=c_2e^{bt}$. By
 (4.32ii) and $b^2-b+\lambda=0$, we get $-b^2+3b=\lambda+2$ and
 $b=1$. But $b<\frac{1}{2}$, this is a contradiction.\\
\noindent {\bf Case ii)(1)(a)2)}~$\underline{b\neq
0},~\underline{p_1\lambda_2=
 0}$\\
If $\underline{p_1=0},$ then $\zeta=2p_2$, $\eta=2p_2^2$ and
$b=\frac{(\lambda+2)\eta}{\zeta^2}$, so $\lambda=-4$ and $b=-1$. By
(4.32iii), we get $\lambda_2=0$ and $-2b^2+4b=\lambda+2$ which is a
contradiction.\\
If $\underline{\lambda_2=0}$, by (4.32ii) and
$b=\frac{(\lambda+2)\eta}{\zeta^2}$, we get $\lambda=0$ or $-2$.
When $\lambda=0$, then $b=\frac{2\eta}{\zeta^2}=0$, this is a
contradiction. There is a similar contradiction for $\lambda=-2$.\\
\noindent {\bf Case ii)(1)(a)3)}~$\underline{b= 0},$\\
 then $\psi=c_2$, by (4.32i), $\lambda=0$. By (4.32ii),
 $\lambda=-2$, this is a contradiction.\\

 \noindent {\bf Case ii)(1)(b)}~$\underline{c_2=0}$,\\
 then $$\left[a-\frac{(\lambda+2)\eta}{\zeta^2}\right]c_1e^{at}=
\frac{p_1\lambda_2\eta}{(p_2-p_1)\zeta^2}(
c_1e^{at})^{1-\frac{2p_2\zeta}{\eta}}.\eqno(4.36)$$ \noindent {\bf
Case ii)(1)(b)1)} $\underline{p_1\lambda_2}\neq 0$,\\
then $p_2=0$ and $\zeta=2p_2$, $\eta=2p_2^2$ and
$a=\frac{1+\sqrt{1-4\lambda}}{2}$ and
$\lambda_2=\lambda+2-\frac{1+\sqrt{1-4\lambda}}{2}.$ By (4.32ii),
then $a=1$ and $\lambda=0$, so $\lambda_2=1$ and $\psi=c_1e^t$ and
$\phi$ satisfies (4.32iii). In this case, we get $\underline{ p_2=0,
~p_1\neq 0,~\phi=c_0e^{\frac{t}{p_1}},~\lambda=0,~\lambda_2=1}.$\\
\noindent {\bf
Case ii)(1)(b)2)} $\underline{p_1\lambda_2}= 0$,\\
if $\underline{p_1=0}$, then $\zeta=2p_2$, $\eta=2p_2^2$ and
$\psi=c_1e^{at}$ and $a=\frac{(\lambda+2)\eta}{\zeta^2}$, so
$\lambda=0$ and $a=1$. By (4.32iii), we get $\lambda_2=0$ and $\phi$
satisfies (4.23ii) and (4.32iii). In this case,
$\underline{p_1=0,~p_2\neq
0~\lambda=0,~\lambda_2=0,\phi=c_0e^{\frac{t}{p_2}}}.$\\
If $\underline{\lambda_2=0}$, by (4.32ii) and
$a=\frac{(\lambda+2)\eta}{\zeta^2}$, then $\lambda=0$ and
$a=\frac{2\eta}{\zeta^2}=1$. By (4.32iii), then $\lambda_2=0$ and
$\phi$ satisfies (4.23ii) and (4.32iii). In this case,
$\underline{p_1\neq0,~p_2\neq 0,~\lambda=\lambda_2=0,}$\\
~$\underline{p_1=4p_2,~\phi=c_0e^{\frac{t}{3p_2}}}.$\\

\noindent {\bf Case ii)(1)(c)}~$\underline{c_1\neq 0,~c_2\neq
0,~b\neq 0},$\\
If $\underline{p_2\neq 0}$, then
$e^{at},e^{bt},(c_1e^{at}+c_2e^{bt})^{1-\frac{2p_2\zeta}{\eta}}$ are
linear independent, by (4.34), then
$$\left[a-\frac{(\lambda+2)\eta}{\zeta^2}\right]c_1=0,~\left[b-\frac{(\lambda+2)\eta}{\zeta^2}\right]c_2=0,
\frac{p_1\lambda_2\eta}{(p_2-p_1)\zeta^2}(
c_1e^{at})^{1-\frac{2p_2\zeta}{\eta}}.\eqno(4.37)$$ So
$a=b=\frac{(\lambda+2)\eta}{\zeta^2}$, this is a contradiction.\\
If $\underline{p_2= 0}$, then by (4.34),
$$a-\frac{(\lambda+2)\eta}{\zeta^2}-\frac{p_1\lambda_2\eta}{(p_2-p_1)\zeta^2}=0,~
b-\frac{(\lambda+2)\eta}{\zeta^2}-\frac{p_1\lambda_2\eta}{(p_2-p_1)\zeta^2}=0,\eqno(4.38)$$
so $a=b$ and we get a contradiction.\\
\noindent {\bf Case ii)(1)(d)}~$\underline{c_1\neq 0,~c_2\neq
0,~b= 0},$\\
When $\underline{1-\frac{2p_2\zeta}{\eta}\neq 0}$, we have similar
discussions. When  $\underline{1-\frac{2p_2\zeta}{\eta}= 0}$, we
have $\frac{(\lambda+2)\eta}{\zeta^2}=1$. By $b=0$, then $\lambda=0$
and $2\eta=\zeta^2=4p_2\zeta$, so $4p_2=\zeta$ and $p_1=2p_2$. But
$\eta=2p_2\zeta$, then $p_1=p_2=0$. This is a contradiction.\\

\noindent {\bf Case ii)(2)}~~$\lambda= \frac{\zeta^2}{4\eta},$
$\psi=c_1e^{\frac{1}{2}t} + c_2te^{\frac{1}{2}t},$\\
by (4.33), we have
$$\left[\frac{1}{2}c_1+c_2-a_0c_1+(\frac{c_2}{2}-a_0c_2)t\right]e^{\frac{1}{2}t}=
\frac{p_1\lambda_2\eta}{(p_2-p_1)\zeta^2}(c_1+c_2t)^{1-\frac{2p_2\zeta}{\eta}}(e^{\frac{1}{2}t})^
{1-\frac{2p_2\zeta}{\eta}},\eqno(4.39)$$ where
$a_0=\frac{1}{4}+\frac{2\eta}{\zeta^2}$.\\
 \noindent {\bf Case ii)(2)a)}~~$\underline{c_2\neq 0},$\\
so $\frac{1}{2}c_1+c_2-a_0c_1+(\frac{c_2}{2}-a_0c_2)t\neq 0$ and
$p_2=0$. By (4.32iii),
$\lambda_2+\frac{(\phi^\zeta)'}{\phi^\zeta}=\lambda+2$, then $
\phi^\zeta=c_0e^{(-\lambda_2+\lambda+2)t}$ and $(c_1e^{\frac{1}{2}t}
+
c_2te^{\frac{1}{2}t})^{\frac{\zeta^2}{\eta}}=c_0e^{(-\lambda_2+\lambda+2)t}$,
this is a contradiction with $c_2\neq 0.$\\
\noindent {\bf Case ii)(2)b)}~~$\underline{c_2= 0},$\\
by (4.39), we have
$$\left(\frac{1}{2}c_1-a_0c_1\right)e^{\frac{1}{2}t}=
\frac{p_1\lambda_2\eta}{(p_2-p_1)\zeta^2}(c_1)^{1-\frac{2p_2\zeta}{\eta}}(e^{\frac{1}{2}t})^
{1-\frac{2p_2\zeta}{\eta}},\eqno(4.40)$$ If
$\underline{a_0=\frac{1}{2}},$ then $p_1\lambda_2=0$,~and
$\lambda=2$. If $\underline{p_1=0}$, then $\zeta=2p_2$ and by
(4.32ii), $\frac{(\phi^\zeta)'}{\phi^\zeta}=4$ and
$\phi^\zeta=c_0e^{4t}$. By (4.32iii),
$$\frac{\lambda_2}{(c_0'e^{\frac{1}{2}t})^{\frac{2p_2\zeta}{\eta}}}
-\frac{16p_2}{\zeta}+4\frac{\zeta+2p_2}{\zeta}=4,\eqno(4.41)$$ so
$\lambda_2=0$ and we have a contradiction by (4.41).\\
If $\underline{p_1\neq 0}$, then $\lambda_2=0$, so
$\psi=c_1e^{\frac{1}{2}t}$ and $\phi^\zeta=c_0e^{4t}. $ Then by
(4.32ii), we have $p_1=0$ which contradicts with $p_1\neq 0$.\\
If $\underline{a_0\neq\frac{1}{2}},$ then $p_2=0$ and
$\zeta=p_1,~\eta=p_1^2$ and $\lambda=\frac{1}{4}$ and
$a_0=\frac{9}{4}$. By $\phi^\xi=c_1e^{\frac{1}{2}t}$ and (4.32ii),
we have a contradiction. In a word, we have no solutions in case
ii)(2).\\

\noindent {\bf Case ii)(3)}~~$\lambda>\frac{\zeta^2}{4\eta},$
$\psi=c_1e^{\frac{1}{2}t}{\rm cos}\left(a t\right)
 + c_2e^{\frac{1}{2}t}{\rm sin}\left(a
t\right),$ where
$a=\frac{\sqrt{\frac{4\lambda\eta}{\zeta^2}-1}}{2}$. By (4.33), we
have
$$(\frac{c_1}{2}+ac_2){\rm cos}(at)+(-ac_1+\frac{c_2}{2}){\rm
sin}(at)=\frac{p_1\lambda_2\eta}{(p_2-p_1)\zeta^2}(c_1{\rm
cos}(at)+c_2{\rm
sin}(at))^{1-\frac{2p_2\zeta}{\eta}}e^{-\frac{2p_2\zeta}{\eta}}$$
$$+ \frac{(\lambda+2)\eta}{\zeta^2}(c_1{\rm cos}(at)+c_2{\rm
sin}(at)).\eqno(4.42)$$ If $\underline{p_2\neq 0},$ then
$p_1\lambda_2=0$ and
$$\frac{c_1}{2}+ac_2=\frac{(\lambda+2)\eta}{\zeta^2}c_1,~-ac_1+\frac{c_2}{2}=\frac{(\lambda+2)\eta}{\zeta^2}c_2,\eqno(4.43)$$
so $c_1^2+c_2^2=0$. This is a contradiction.\\
If $\underline{p_2= 0},$ then
$$\frac{c_1}{2}+ac_2=\frac{p_1\lambda_2\eta}{(p_2-p_1)\zeta^2}c_1+\frac{(\lambda+2)\eta}{\zeta^2}c_1,~-ac_1+\frac{c_2}{2}
=\frac{p_1\lambda_2\eta}{(p_2-p_1)\zeta^2}c_2+\frac{(\lambda+2)\eta}{\zeta^2}c_2,\eqno(4.44)$$
Then $c_1^2+c_2^2=0$. This is a contradiction. By the above discussions, we get the following theorem:\\

\noindent {\bf Theorem 4.28}{\it~Let
$M=I\times_{\phi^{p_1}}F_1\times_{\phi^{p_2}}F_2$ be a generalized
Kasner spacetime and ${\rm dim}F_1=1,~{\rm dim}F_2=2$ and
$P=\frac{\partial}{\partial t}$. Then
    $(M,\overline{\nabla})$ is Einstein with the Einstein constant
    $\lambda$ if and only if
 $(F_2,\nabla^{F_2})$ is Einstein with the Einstein
constant $\lambda_2$, and  one of the following conditions is satisfied}\\
\noindent (1) $ p_2=0, ~p_1\neq
0,~\phi=c_0e^{\frac{t}{p_1}},~\lambda=0,~\lambda_2=1.$\\
\noindent (2) $p_1=0,~p_2\neq
0~\lambda=0,~\lambda_2=0,\phi=c_0e^{\frac{t}{p_2}}.$\\
\noindent (3) $p_1\neq0,~p_2\neq 0,~\lambda=\lambda_2=0,~
,p_1=4p_2,~\phi=c_0e^{\frac{t}{3p_2}}.$\\

\noindent{\bf $\cdot$ Type (II) generalized Kasner space-times with
a semi-symmetric metric
 connection with constant scalar curvature}\\
\indent By Proposition 4.26, then $(F_2,\nabla^{F_2})$ has constant
    scalar curvature $S^{F_2}$ and
$$\overline{S}=\frac{S^{F_2}}{\phi^{2p_2}}-2\zeta\frac{\phi''}{\phi}-(\eta+\zeta^2-2\zeta)\frac{\phi'^2}{\phi^2}+6\zeta
\frac{\phi'}{\phi}-6.\eqno(4.45)$$ If $\underline{\zeta=0}$, when
$\underline{\eta=0},$  then $p_1=p_2=0$ and
$\underline{\overline{S}=S^{F_2}-6}$. If $\underline{\eta\neq 0},$
then
$$\eta\frac{\phi'^2}{\phi^2}=\frac{S^{F_2}}{\phi^{2p_2}}-(\overline{S}+6).\eqno(4.46)$$
If $\underline{\zeta\neq 0},$ putting
$\phi=\psi^{\frac{2\zeta}{\eta+\zeta^2}},$ we get
$$-\frac{4\zeta^2}{\eta+\zeta^2}\psi''+\frac{12\zeta^2}{\eta+\zeta^2}\psi'-(\overline{S}+6)\psi
+S^{F_2}\psi^{1-\frac{4p_2\zeta}{\eta+\zeta^2}}=0.\eqno(4.47)$$

\noindent{\bf $\cdot$ Type (III) generalized Kasner space-times with
a semi-symmetric metric
 connection with constant scalar curvature}\\

\indent By Proposition 4.26, then
$$\overline{S}=-2\zeta\frac{\phi''}{\phi}-(\eta+\zeta^2-2\zeta)\frac{\phi'^2}{\phi^2}+6\zeta
\frac{\phi'}{\phi}-6.\eqno(4.48)$$ If $\underline{\zeta=\eta=0},$
then $p_1=p_2=p_3=0$, we get $\overline{S}=-6$.\\
If $\underline{\zeta=0,~\eta\neq 0},$ then $[({\rm
ln}\phi)']^2=-\frac{\overline{S}+6}{\eta}$, so when
$\overline{S}+6>0$, there is no solutions, when $\overline{S}+6=0$,
$\phi$ is a constant and  when $\overline{S}+6<0$,
$\phi=c_0e^{\pm\sqrt{-\frac{\overline{S}+6}{\eta}}t}.$\\
If $\underline{\zeta\neq 0}$, then $\eta\neq 0$, putting
$\phi=\psi^{\frac{2\zeta}{\eta+\zeta^2}}$, then
$$\psi''-3\psi'+\frac{(\overline{S}+6)(\eta+\zeta^2)}{4\zeta^2}\psi=0.\eqno(4.49)$$
So, we get \\
\noindent (1) $\overline{S}+6<\frac{9\zeta^2}{\eta+\zeta^2},$~
$\psi=c_1e^{\frac{3+\sqrt{9-\frac{(\overline{S}+6)(\eta+\zeta^2)}{\zeta^2}}}{2}t}+c_2
e^{\frac{3-\sqrt{9-\frac{(\overline{S}+6)(\eta+\zeta^2)}{\zeta^2}}}{2}t},$\\
\noindent (2) $\overline{S}+6=\frac{9\zeta^2}{\eta+\zeta^2},$~
$\psi=c_1e^{\frac{3}{2}t}+c_2te^{\frac{3}{2}t},$\\
\noindent (3) $\overline{S}+6>\frac{9\zeta^2}{\eta+\zeta^2},$~
$\psi=c_1e^{\frac{3}{2}t}{\rm
cos}\left(\frac{\sqrt{-9+\frac{(\overline{S}+6)(\eta+\zeta^2)}{\zeta^2}}}{2}t\right)+c_2e^{\frac{3}{2}t}
{\rm
sin}\left(\frac{\sqrt{-9+\frac{(\overline{S}+6)(\eta+\zeta^2)}{\zeta^2}}}{2}t\right).$
So we get the following theorem\\

\noindent {\bf Theorem 4.29}{\it~Let
$M=I\times_{\phi^{p_1}}F_1\times_{\phi^{p_2}}F_2\times_{\phi^{p_3}}F_3$
be a generalized Kasner spacetime and ${\rm dim}F_1={\rm
dim}F_2={\rm dim}F_3=1,$ and $P=\frac{\partial}{\partial t}$. Then
$\overline{S}$ is a constant if and only if one of the following
case holds} \\
\noindent (1) $\zeta=\eta=0,$ $\overline{S}=-6$.\\
\noindent (2){\it $ \zeta=0,~\eta\neq 0,$ when $\overline{S}+6>0$,
there is no solutions, when $\overline{S}+6=0$, $\phi$ is a constant
and
when $\overline{S}+6<0$, $\phi=c_0e^{\pm\sqrt{-\frac{\overline{S}+6}{\eta}}t}.$}\\
\noindent (3) If $\zeta\neq 0$\\
 \noindent (3a)
$\overline{S}+6<\frac{9\zeta^2}{\eta+\zeta^2},$~
$\phi=\left(c_1e^{\frac{3+\sqrt{9-\frac{(\overline{S}+6)(\eta+\zeta^2)}{\zeta^2}}}{2}t}+c_2
e^{\frac{3-\sqrt{9-\frac{(\overline{S}+6)(\eta+\zeta^2)}{\zeta^2}}}{2}t}\right)^{\frac{2\zeta}{\eta+\zeta^2}},$\\
\noindent (3b) $\overline{S}+6=\frac{9\zeta^2}{\eta+\zeta^2},$~
$\phi=\left(c_1e^{\frac{3}{2}t}+c_2te^{\frac{3}{2}t}\right)^{\frac{2\zeta}{\eta+\zeta^2}},$\\
\noindent (3c) $\overline{S}+6>\frac{9\zeta^2}{\eta+\zeta^2},$~
$\phi=\left(c_1e^{\frac{3}{2}t}{\rm
cos}\left(\frac{\sqrt{-9+\frac{(\overline{S}+6)(\eta+\zeta^2)}{\zeta^2}}}{2}t\right)+c_2e^{\frac{3}{2}t}
{\rm
sin}\left(\frac{\sqrt{-9+\frac{(\overline{S}+6)(\eta+\zeta^2)}{\zeta^2}}}{2}t\right)\right)^{\frac{2\zeta}{\eta+\zeta^2}}.$\\

\noindent{\bf $\cdot$ Einstein Type (III) generalized Kasner
space-times with a semi-symmetric metric
 connection}\\

\indent By Proposition 4.25, we have\\

$$\zeta\left(\frac{\phi'-\phi''}{\phi}\right)-(\eta-\zeta)\frac{\phi'^2}{\phi^2}=\lambda,\eqno(4.50i)$$
$$-p_1\left[\frac{\phi''}{\phi}+(\zeta-1)\frac{\phi'^2}{\phi^2}-2\frac{\phi'}{\phi}\right]+\zeta\frac{\phi'}{\phi}=\lambda+2,\eqno(4.50ii)$$
$$-p_2\left[\frac{\phi''}{\phi}+(\zeta-1)\frac{\phi'^2}{\phi^2}-2\frac{\phi'}{\phi}\right]+\zeta\frac{\phi'}{\phi}=\lambda+2,\eqno(4.50iii)$$
$$-p_3\left[\frac{\phi''}{\phi}+(\zeta-1)\frac{\phi'^2}{\phi^2}-2\frac{\phi'}{\phi}\right]+\zeta\frac{\phi'}{\phi}=\lambda+2,\eqno(4.50iv)$$
If $\underline{\zeta=\eta=0},$ by (4.50i), $\lambda=0$, by (4.50ii),
$\lambda=-2$, this is a contradiction.\\
If $\underline{\zeta=0, ~\eta\neq0},$ plusing
(4.50ii),(4.50iii),(4.50iv), we get $\lambda=-2$. By (4.50i),
$\frac{\phi'^2}{\phi^2}=\frac{2}{\eta}$ and
$\phi=c_0e^{\pm\sqrt{\frac{2}{\eta}}t}.$ But by (4.50ii), then
$\frac{\phi''}{\phi}+(\zeta-1)\frac{\phi'^2}{\phi^2}-2\frac{\phi'}{\phi}=0$,
this is a contradiction.\\
$\underline{\zeta\neq 0}$. If $p_1=p_2=p_3$, we get type (I), so we
may let $p_1\neq p_2$. By (4.50ii) and (4.50iii), we have
$\frac{(\phi^\zeta)'}{\phi^\zeta}=\lambda+2$ and
$\frac{(\phi^\zeta)''}{\phi^\zeta}-2\frac{(\phi^\zeta)'}{\phi^\zeta}=0,$
so $\phi^\zeta=c_0e^{(\lambda+2)t}$ and $\lambda=-2$ or $0$. When
$\lambda=-2$, $\psi$ is a constant, by (4.50i), $\lambda=0$, this is
a contradiction. When $\lambda=0$, $\psi=c_0e^{\frac{2\eta
t}{\zeta^2}}$ and $\psi''-\psi'=0$, so $\frac{2\eta }{\zeta^2}=1$.
In this case, we get when $p_i\neq p_j$ for some $i,j\in\{1,2,3\}$,
$\underline{\lambda=0,~\frac{2\eta
}{\zeta^2}=1,~\phi=c_0e^{\frac{2t}{\zeta}},~}.$ We get the following
theorem.\\

\noindent {\bf Theorem 4.30}{\it~Let
$M=I\times_{\phi^{p_1}}F_1\times_{\phi^{p_2}}F_2\times_{\phi^{p_3}}F_3$
be a generalized Kasner spacetime for $p_i\neq p_j$ for some
$i,j\in\{1,2,3\}$ and ${\rm dim}F_1={\rm dim}F_2={\rm dim}F_3=1,$
and $P=\frac{\partial}{\partial t}$.
 Then
    $(M,\overline{\nabla})$ is Einstein with the Einstein constant
    $\lambda$ if and only if $\lambda=0,~\frac{2\eta
}{\zeta^2}=1,~\phi=c_0e^{\frac{2t}{\zeta}}.$}\\

\section{Multiply twisted product Finsler
 manifolds}

 \quad In this section, we set $(M_i,F_i)$ is a Finsler manifold
 for $0\leq i\leq {\bf b}$ and $f_i:M_0\times M_i\rightarrow
 \mathbb{R}$ is a smooth function for $1\leq i\leq {\bf b}$
Let $\pi _i: TM_i\rightarrow M_i$ be the projection map. The product
manifold $M_0\times M_1\times \cdots \times M_{\bf b}$
      endowed with the metric ${F}:TM^0\times TM_1^0\cdots\times
      TM_{\bf b}^0\rightarrow \mathbb{R}$ is considered,\\
      $~~~~~{F}(v_0,v_1,\cdots,v_m)$\\
      $$=\sqrt{{F_0}^2(v_0)+f_1^2(\pi_0(v_0),\pi_1(v_1)){F_1}^2(v_1)+\cdots+
     f_m^2(\pi_0(v_0),\pi_m(v_m)){F_{\bf b}}^2(v_{\bf b})},\eqno(5.1)$$
     where $TM^0_i=TM_i-\{0\}.$ Let ${\rm dim M_i}=m_i,$ for $0\leq i\leq {\bf b}$
and $(x_i^1,\cdots,x_i^{m_i},y_i^1,\cdots,y_i^{m_i}$ is the local
coordinate on $TM_i$. Let
$$g_{ij}=\frac{1}{2}\frac{\partial^2F^2}{\partial y^i\partial y^j},\eqno(5.2)$$
and $g^{ij}$ be the inverse of $g_{ij}$. For a Finsler manifold
$(M,F)$, a global vector field in introduced by $F$ on $TM^0$, which
in a standard coordinate $(x^i,y^i)$ for $TM^0$ is given by ${\bf
G}=y^i\frac{\partial}{\partial
x^i}-2G^i(x,y)\frac{\partial}{\partial y^i}$, where
$$G^i:=\frac{1}{4}g^{il}\left(\frac{\partial^2F^2}{\partial x^k\partial
y^l}y^k-\frac{\partial F^2}{\partial x^l}\right).\eqno(5.3)$$ Let
$$G^i_j=\frac{\partial G^i}{\partial y^j};~
G^i_{jk}=\frac{\partial^2 G^i}{\partial y^j\partial
y^k};~B^{i}_{jkl}=\frac{\partial^3 G^i}{\partial y^j\partial
y^k\partial y^l};~E_{jk}=\frac{1}{2}B^{i}_{jki}.\eqno(5.4)$$ $F$ is
called a Berwald metric and weakly Berwald metric if $B^{i}_{jkl}=0$
$E_{jk}=0$ respectively. Let $C_{ijk}=\frac{1}{2}\frac{\partial
g_{ij}}{\partial y^k}.$ A Finsler metric $F$ is said to be isotropic
mean Berwald metric if its mean Berwald curvature is in the
following form
$$E_{ij}=\frac{1}{2}(n+1)cF^{-1}h_{ij},\eqno(5.5)$$
where $h_{ij}=g_{ij}-F^{-2}y_iy_j$ is the angular metric and
$c=c(x)$ is a scalar function on $M$. Let
$L_{ijk}=\frac{-1}{2}y_lB^l_{ijk}$. A Finsler metric is called a
Landsgerg metric if $L_{ijk}=0$. A Finsler metric is said to be
relatively isotropic Landsberg metric if it satisfies
$L_{ijk}=cFC_{ijk}$, where $c=c(x)$ is a scalar function on $M$. Let
$J_i=g^{jk}L_{ijk}$ and $I_i=g^{jk}C_{ijk}.$ A Finsler metric is
called a weakly Landsberg metric if $J_i=0$. A Finsler metric is
said to be relatively isotropic mean Landsberg metric if $J_i=cFI_i$
for some scalar function $c=c(x)$ on $M$. A Finsler metric is said
to be locally dually flat if $(F^2)_{x^ky^l}y^k=2(F^2)_{x^l}.$ A
Finsler manifold $(M,F)$ is called a locally Minkowski manifold if
on each coordinate neighborhood of $TM$, $F$ is a function of
$(y^i)$ only. Let $I_i=g^{jk}C_{ijk}$ and
$h_{ij}=g_{ij}-\frac{1}{F^2}g_{ip}y^pg_{iq}y^q$ and
$M_{ijk}=C_{ijk}-\frac{1}{n+1}(I_ih_{jk}+I_jh_{ik}+I_kh_{ij})$. A
Finsler metric $F$ is said to be $C$-reducible if $M_{ijk}=0.$
 Direct computations show that
$$G^j=(G_0)^j-\frac{1}{4}\sum_{i=1}^{{\bf
b}}\sum_{k=1}^{m_0}g_0^{jk}\frac{\partial f_i^2}{\partial
x_0^k}F_i^2,~{\rm for} ~1\leq j\leq m_0,\eqno(5.6)$$ where $(G_l)^j$
is the geodesic coefficient of $(M_l,F_l)$ and
$(g_l)^{jk}=\frac{1}{2}\frac{\partial^2F_l^2}{\partial y^i_l\partial
y^j_l}$ for $0\leq l\leq {\bf b}.$
$$G^{m_0+\cdots+m_{l-1}+j}=(G_l)^j+\frac{1}{2f_l^2}y_0^hy_l^j\frac{\partial
f_l^2}{\partial x_0^h}+\frac{1}{2f_l^2}y_l^ry_l^j\frac{\partial
f_l^2}{\partial x_l^r}-\frac{1}{4f_l^2}(g_l)^{js}\frac{\partial
f_l^2}{\partial x_l^s}F_l^2,\eqno(5.7)$$ for $1\leq l\leq {\bf
b},~1\leq j\leq m_l.$ In the following, we write
$\overline{l-1}=m_0+\cdots+m_{l-1}.$
$$G^j_k=(G_0)_k^j-\frac{1}{4}\sum_{i=1}^{{\bf
b}}\sum_{k=1}^{m_0}\frac{\partial g_0^{js}}{\partial
y_0^k}\frac{\partial f_i^2}{\partial x_0^s}F_i^2,\eqno(5.8)$$
$$G^j_{\overline{r-1}+k}=-\frac{1}{2}\sum_{s=1}^{m_0}g_0^{js}\frac{\partial
f_r^2}{\partial x_0^s}(g_r)_{kt}y_r^t,\eqno(5.9)$$
$$G_k^{\overline{l-1}+j}=\frac{1}{2f_l^2}y_l^j\frac{\partial
f_l^2}{\partial x_0^k},\eqno(5.10)$$
$$G_{\overline{r-1}+k}^{\overline{l-1}+j}=\delta^l_r(G_l)^j_k+\frac{1}{2f_l^2}y_0^h
\frac{\partial f_l^2}{\partial
x_0^h}\delta^k_j\delta^l_r+\delta^l_rL_k^j,\eqno(5.11)$$ where
$$L_k^j=
\frac{1}{2f_l^2}\delta^k_jy^r_l\frac{\partial f_l^2}{\partial x_l^r}
+\frac{1}{2f_l^2}y^j_l\frac{\partial f_l^2}{\partial
x_l^k}-\frac{1}{4f_l^2}\frac{\partial f_l^2}{\partial
x_l^s}\left(\frac{\partial (g_l)^{js}}{\partial
y_l^k}F_l^2+2(g_l)^{js}(g_l)_{kt}y_l^t\right).$$

$$G^k_{ij}=(G_0)_{ij}^k-\frac{1}{4}\sum_{t=1}^{{\bf
b}}\sum_{k=1}^{m_0}\frac{\partial^2 g_0^{ks}}{\partial y_0^i\partial
y_0^j}\frac{\partial f_t^2}{\partial x_0^s}F_t^2,\eqno(5.12)$$

$$G^k_{\overline{r-1}+i,j}=-\frac{1}{2}\sum_{s=1}^{m_0}\frac{\partial g_0^{ks}}{\partial y_0^j}\frac{\partial
f_r^2}{\partial x_0^s}(g_r)_{it}y_r^t,\eqno(5.13)$$

$$G_{ij}^{\overline{r-1}+k}=0,~G_{\overline{l-1}+i,j}^{\overline{r-1}+k}=\frac{1}{2f_r^2}\delta_i^k\delta_r^l
\frac{\partial f_r^2}{\partial x_0^j},\eqno(5.14)$$
$$G_{\overline{l-1}+i,\overline{r-1}+j}^k=-\frac{1}{2}(g_0)^{ks}\frac{\partial
f_r^2}{\partial x_0^s}\delta^l_r(g_r)_{ij},\eqno(5.15)$$
$$G_{\overline{s-1}+i,\overline{r-1}+j}^{\overline{l-1}+k}=\delta_r^l\delta_s^l(G_l)^k_{ij}+\frac{1}{2f_l^2}\delta_r^l\delta_j^k\delta_s^l
\frac{\partial f_r^2}{\partial
x_l^i}+\frac{1}{2f_l^2}\delta_s^l\delta_i^k\delta_r^l\frac{\partial
f_l^2}{\partial x_l^j}-\frac{1}{4f_l^2}\frac{\partial
f_l^2}{\partial
x_l^\beta}\delta_r^l\delta_s^l\frac{\partial^2[(g_l)^{k\beta}F_l^2]}{\partial
y_l^i \partial y_l^j}.\eqno(5.16)$$
$$C_{ijk}=(C_0)_{ijk},~C_{\overline{s-1}+i,\overline{r-1}+j,
\overline{l-1}+k}=\delta_r^s\delta^l_rf_s^2(C_r)_{ijk},\eqno(5.17)$$
and the other six combinations are zero. Let
$C_{ij}^k=\frac{1}{2}g^{kh}\frac{\partial g_{ij}}{\partial y^h}.$
Then
$$C_{ij}^k=(C_0)_{ij}^k,~C_{\overline{s-1}+i,\overline{r-1}+j}
^{\overline{l-1}+k}=\delta_r^s\delta^l_rf_s^2(C_l)_{ij}^k,\eqno(5.18)$$
and the other six combinations are zero. Let
$$F_{ij}^k=\frac{1}{2}g^{kh}\left(\frac{\delta g_{hi}}{\delta x^j}+
\frac{\delta g_{hj}}{\delta x^i}-\frac{\delta g_{ij}}{\delta
x^h}\right),\eqno(5.19)$$ where
$$\frac{\delta}{\delta
x^k}=\frac{\partial}{\partial x_0^k}-\sum_{l=0}^{{\bf b}}\sum
_{j_l=1}^{m_l}G_k^{\overline{l-1}+j_l}\frac{\partial}{\partial
y_l^{j_l}},~1\leq k \leq m_0,\eqno(5.20)$$
$$\frac{\delta}{\delta
x^{\overline{r-1}+k}}=\frac{\partial}{\partial x_r^k}-
\sum_{l=0}^{{\bf b}}\sum
_{j_l=1}^{m_l}G_{\overline{r-1}+k}^{\overline{l-1}+j_l}\frac{\partial}{\partial
y_l^{j_l}},~1\leq k \leq m_r,\eqno(5.21)$$ Then
$$F_{ij}^k=(F_0)_{ij}^k+\frac{1}{4}(g_0)^{kh}\left[\frac{\partial
g_0^{ls}}{\partial y_0^j}(C_0)_{hil}+\frac{\partial
g_0^{ls}}{\partial y_0^i}(C_0)_{hjl}-\frac{\partial
g_0^{ls}}{\partial y_0^h}(C_0)_{jil}\right](\sum _{i=1}^{{\bf
b}}\frac{\partial f_r^2}{\partial x_0^s}F_r^2),\eqno(5.22)$$
$$F_{ij}^{\overline{l-1}+k}=-\frac{1}{4}f_l^{-2}(g_0)^{s\beta}\frac{\partial
f_l^2}{\partial x_0^\beta}y_l^k\frac{\partial(g_0)_{ij}}{\partial
y_0^s},\eqno(5.23)$$
$$F_{i,\overline{r-1}+j}^k=-\frac{1}{4}\frac{\partial(g_0)^{k\beta}}{\partial
y_0^i}\frac{\partial f_r^2}{\partial
x_0^\beta}y_r^t(g_r)_{tj},\eqno(5.24)$$
$$F_{i,\overline{r-1}+j}^{\overline{l-1}+k}=\frac{1}{2f_l^2}\delta_j^k\delta_r^l
\frac{\partial f_l^2}{\partial x_0^i},\eqno(5.25)$$
$$F_{\overline{r-1}+i,j}^k=F_{j,\overline{r-1}+i}^k,~~
F_{\overline{r-1}+j,i}^{\overline{l-1}+k}=F_{i,\overline{r-1}+j}^{\overline{l-1}+k},\eqno(5.26)$$
$$F_{\overline{s-1}+i,\overline{r-1}+j}^k=-\frac{1}{2}\frac{\partial f_r^2}{\partial
x_0^h}g_0^{kh}(g_r)_{ij}\delta_r^s,\eqno(5.27)$$

$$F_{\overline{s-1}+i,\overline{r-1}+j}^{\overline{l-1}+k}=
\frac{1}{2f_l^2}\delta_s^l\delta_r^l\delta_k^i\frac{\partial
f_l^2}{\partial x_l^j}+
\frac{1}{2f_l^2}\delta_s^l\delta_r^l\delta_k^j\frac{\partial
f_l^2}{\partial x_l^i}-
\frac{1}{2f_l^2}\delta_s^l\delta_r^l(g_l)^{kh}(g_l)_{ij}
\frac{\partial f_l^2}{\partial
x_l^h}+\delta_s^l\delta_r^l(F_l)_{ij}^k$$
$$-\frac{1}{4f_l^2}(g_l)^{kh}\delta_s^l\delta_r^l \frac{\partial
f_l^2}{\partial x_0^d}y_0^d\frac{\partial (g_l)_{ij}}{\partial
y_l^h}-\frac{1}{2}(g_l)^{kh}\delta_s^l\delta_r^l\left(L_j^\beta
\frac{\partial (g_l)_{hi}}{\partial y_l^\beta}+ L_i^\beta
\frac{\partial (g_l)_{hj}}{\partial y_l^\beta}-L_h^\beta
\frac{\partial (g_l)_{ij}}{\partial y_l^\beta}\right).\eqno(5.28)$$

$$B^k_{ijq}=(B_0)_{ijq}^k-\frac{1}{4}\sum_{t=1}^{{\bf
b}}\sum_{k=1}^{m_0}\frac{\partial^3 g_0^{ks}}{\partial y_0^i\partial
y_0^j\partial y_0^q}\frac{\partial f_t^2}{\partial
x_0^s}F_t^2,\eqno(5.29)$$

$$B^k_{\overline{r-1}+i,j,q}=-\frac{1}{2}\sum_{s=1}^{m_0}\frac{\partial^2 g_0^{ks}}{\partial y_0^j\partial y_0^q}\frac{\partial
f_r^2}{\partial
x_0^s}(g_r)_{it}y_r^t=B^k_{j,\overline{r-1}+i,q}=B^k_{j,q,\overline{r-1}+i},\eqno(5.30)$$

$$B_{\overline{l-1}+i,\overline{r-1}+j,q}^k=-\frac{1}{2}\frac{\partial (g_0)^{ks}}{\partial y_0^q}\frac{\partial
f_r^2}{\partial x_0^s}\delta^l_r(g_r)_{ij},\eqno(5.31)$$

$$B_{\overline{l-1}+i,\overline{r-1}+j,\overline{s-1}+q}^k=- (g_0)^{ks}\frac{\partial
f_r^2}{\partial x_0^s}\delta^l_r\delta^l_s(C_r)_{ijq},\eqno(5.32)$$

$$B_{i,j,q}^{\overline{r-1}+k}=B_{\overline{l-1}+i,j,q}^{\overline{r-1}+k}
=B_{\overline{l-1}+i,\overline{s-1}+j,q}^{\overline{r-1}+k}
=0,\eqno(5.33)$$

$$B_{\overline{s-1}+i,\overline{r-1}+j,\overline{t-1}+q}^{\overline{l-1}+k}=\delta_r^l\delta_s^l\delta_t^l
(B_l)^k_{ijq}-\frac{1}{4f_l^2}\frac{\partial f_l^2}{\partial
x_l^\beta}\delta_r^l\delta_s^l\delta_t^l\frac{\partial^3[(g_l)^{k\beta}F_l^2]}{\partial
y_l^i \partial y_l^j\partial y_l^q}.\eqno(5.34)$$ By (5.18), then\\

\noindent{\bf Proposition 5.1} {\it The multiply twisted product
Finsler manifold $M_0\times_{f_1} M_1\times \cdots \times_{f_{\bf
b}} M_{\bf b}$ is Riemannian if and only if $(M_i,F_i)$ is
Riemannian for $0\leq i
\leq {\bf b}.$}\\

Similar to Theorem 1 in [PTN], by (5.29)-(5.34) we have\\

\noindent{\bf Theorem 5.2} {\it If there is a point $p_i\in M_0$
such that $\frac{\partial f_i^2}{\partial x_0^j(p_i,x_i)}\neq 0$ for
some $j\in \{1,\cdots,m_0\}$, then $M_0\times_{f_1} M_1\times \cdots
\times_{f_{\bf b}} M_{\bf b}$ is a Berwald manifold if and only if
$(M_0,F_0)$ is a Berwald manifold, $(M_i,F_i)$ is Riemannian for
$1\leq i \leq {\bf b}$ and $\frac{\partial (g_0)^{ks}}{\partial
y_0^q}\frac{\partial f_l^2}{\partial x_0^s}=0,~1\leq l \leq {\bf
b}.$}\\

By (5.29)-(5.34) then
$$E_{ij}=(E_0)_{ij}-\frac{1}{8}\sum_{t=1}^{{\bf
b}}\sum_{k=1}^{m_0}\frac{\partial^3 g_0^{ks}}{\partial y_0^i\partial
y_0^j\partial y_0^k}\frac{\partial f_t^2}{\partial
x_0^s}F_t^2,\eqno(5.35)$$

$$E_{\overline{l-1}+i,j}=E_{j,\overline{l-1}+i}=
-\frac{1}{4}\sum_{s=1}^{m_0}\frac{\partial^2 g_0^{ks}}{\partial
y_0^j\partial y_0^k}\frac{\partial f_l^2}{\partial
x_0^s}(y_l)_i,\eqno(5.36)$$

$$E_{\overline{l-1}+i,\overline{r-1}+j}=-\frac{1}{4}\frac{\partial (g_0)^{ks}}{\partial y_0^k}\frac{\partial
f_r^2}{\partial x_0^s}\delta^l_r(g_r)_{ij}+\delta^l_r(E_l)_{ij}
-\frac{1}{8f_l^2}\frac{\partial f_l^2}{\partial
x_l^\beta}\delta_r^l\frac{\partial^3[(g_l)^{k\beta}F_l^2]}{\partial
y_l^i \partial y_l^j\partial y_l^k}.\eqno(5.37)$$ By (5.35)-(5.37),
similar to Theorem in [PTN], we have\\

\noindent{\bf Theorem 5.3} {\it $M_0\times_{f_1} M_1\times \cdots
\times_{f_{\bf b}} M_{\bf b}$ is a weakly Berwald manifold if and
only if $(M_0,F_0)$ is a weakly Berwald manifold and $\frac{\partial
(g_0)^{ks}}{\partial y_0^k}\frac{\partial f_l^2}{\partial x_0^s}=0$
and $$(E_l)_{ij}=\frac{1}{8f_l^2}\frac{\partial f_l^2}{\partial
x_l^\beta}\frac{\partial^3[(g_l)^{k\beta}F_l^2]}{\partial y_l^i
\partial y_l^j\partial y_l^k},\eqno(5.38)$$
for $~1\leq l \leq {\bf
b}.$}\\

By (5.5) and (5.35)-(5.37), then\\

\noindent{\bf Lemma 5.4} {\it  $M_0\times_{f_1} M_1\times \cdots
\times_{f_{\bf b}} M_{\bf b}$ has isotropic mean Berwald curvature
if and only if } $$(E_0)_{ij}-\frac{1}{8}\sum_{t=1}^{{\bf
b}}\sum_{k=1}^{m_0}\frac{\partial^3 g_0^{ks}}{\partial y_0^i\partial
y_0^j\partial y_0^k}\frac{\partial f_t^2}{\partial
x_0^s}F_t^2=\frac{1}{2}(n+1)cF^{-1}[(g_0)_{ij}-F^{-2}(y_0)_i(y_0)_j],\eqno(5.39)$$

$$-\frac{1}{4}\sum_{s=1}^{m_0}\frac{\partial^2 g_0^{ks}}{\partial
y_0^j\partial y_0^k}\frac{\partial f_l^2}{\partial
x_0^s}(y_l)_i=-\frac{n+1}{2}cF^{-3}f_l^2(y_l)_i(y_0)_j,\eqno(5.40)$$

$$-\frac{1}{4}\frac{\partial (g_0)^{ks}}{\partial
y_0^k}\frac{\partial f_r^2}{\partial
x_0^s}\delta^l_r(g_r)_{ij}+\delta^l_r(E_l)_{ij}
-\frac{1}{8f_l^2}\frac{\partial f_l^2}{\partial
x_l^\beta}\delta_r^l\frac{\partial^3[(g_l)^{k\beta}F_l^2]}{\partial
y_l^i \partial y_l^j\partial y_l^k}$$
$$=\frac{n+1}{2}cF^{-1}[f_l^2\delta_r^l(g_l)_{ij}-F^{-2}f_l^2f_r^2(y_l)_i(y_r)_j]
.\eqno(5.41)$$

By Lemma 5.4, similar to Theorem 3 in [PTN], we get\\

\noindent{\bf Theorem 5.5} {\it  $M_0\times_{f_1} M_1\times \cdots
\times_{f_{\bf b}} M_{\bf b}$ with isotropic mean Berwald curvature
is a weakly Berwald manifold.}\\

By (5.29)-(5.34), we have
$$L_{ijq}=(L_0)_{ijq}+\frac{1}{8}\sum_{t=1}^{{\bf
b}}\sum_{k=1}^{m_0}\frac{\partial^3 g_0^{ks}}{\partial y_0^i\partial
y_0^j\partial y_0^q}\frac{\partial f_t^2}{\partial
x_0^s}F_t^2(y_0)_k,\eqno(5.42)$$

$$L_{\overline{r-1}+i,j,q}=\frac{1}{4}\sum_{s=1}^{m_0}\frac{\partial^2 g_0^{ks}}{\partial y_0^j\partial y_0^q}\frac{\partial
f_r^2}{\partial x_0^s}(y_0)_k(y_r)_i,\eqno(5.43)$$

$$L_{\overline{l-1}+i,\overline{r-1}+j,q}=0,\eqno(5.44)$$
$$L_{\overline{l-1}+i,\overline{r-1}+j,\overline{t-1}+q}=
\frac{1}{2}y_0^s\frac{\partial f_r^2}{\partial
x_0^s}\delta_r^l\delta_t^l(C_r)_{ijq}+f_l^2\delta_r^l\delta_t^l(L_l)_{ijq}
+\frac{1}{8}(y_l)_k\frac{\partial f_l^2}{\partial
x_l^\beta}\delta_r^l\delta_t^l\frac{\partial^3[(g_l)^{k\beta}F_l^2]}{\partial
y_l^i \partial y_l^j\partial y_l^q}.\eqno(5.45)$$

By (5.42)-(5.45), similar to Theorem 6 in [PTN], we have\\

\noindent{\bf Theorem 5.6} {\it If there is a point $p_i\in M_0$
such that $\frac{\partial f_i^2}{\partial x_0^j(p_i,x_i)}\neq 0$ for
some $j\in \{1,\cdots,m_0\}$, then $M_0\times_{f_1} M_1\times \cdots
\times_{f_{\bf b}} M_{\bf b}$ is a Landsberg manifold if and only if
$(M_0,F_0)$ is Landsberg, $M_i$ is Riemannian for $1\leq i \leq {\bf
b}$ and $\frac{\partial^3 g_0^{ks}}{\partial y_0^i\partial
y_0^j\partial y_0^q}\frac{\partial f_r^2}{\partial
x_0^s}(y_0)_k=0.$}

Similar to Theorem 6 in [PTN], we have\\

\noindent{\bf Theorem 5.7} {\it $M_0\times_{f_1} M_1\times \cdots
\times_{f_{\bf b}} M_{\bf b}$ is a relatively isotropic Landsberg
manifold, then $M$ is a Landsberg manifold.}\\

By (5.42)-(5.45), then
$$J_{i}=(J_0)_{i}+\frac{1}{8}\sum_{t=1}^{{\bf
b}}\sum_{k=1}^{m_0}\frac{\partial^3 g_0^{ks}}{\partial y_0^i\partial
y_0^j\partial y_0^q}\frac{\partial f_t^2}{\partial
x_0^s}F_t^2(g_0)^{jq}(y_0)_k,\eqno(5.46)$$

$$J_{\overline{l-1}+i}=(J_l)_i+\frac{1}{4}\sum_{s=1}^{m_0}\frac{\partial^2 g_0^{ks}}{\partial y_0^j\partial y_0^q}\frac{\partial
f_l^2}{\partial x_0^s}(y_0)_k(y_l)_ig_0^{jq}$$
$$
+\frac{1}{2f_l^2}y_0^s\frac{\partial f_l^2}{\partial
x_0^s}(g_l)^{jq}(C_l)_{ijq}
+\frac{1}{8f_l^2}(g_l)^{jq}(y_l)_k\frac{\partial f_l^2}{\partial
x_l^\beta}\frac{\partial^3[(g_l)^{k\beta}F_l^2]}{\partial y_l^i
\partial y_l^j\partial y_l^q}.\eqno(5.47)$$

Similar to Theorem 8 in [PTN], by (5.46) and (5.47) we have\\

\noindent{\bf Theorem 5.8} {\it If there is a point $p_i\in M_0$
such that $\frac{\partial f_i^2}{\partial x_0^j(p_i,x_i)}\neq 0$ for
some $j\in \{1,\cdots,m_0\}$, then $M_0\times_{f_1} M_1\times \cdots
\times_{f_{\bf b}} M_{\bf b}$ is a weakly Landsberg manifold if and
only if $(M_0,F_0)$ is a weakly Landsberg manifold, $(M_i,F_i)$ is
Riemannian for $1\leq i \leq {\bf b}$ and
$(y_0)_k(g_0)^{jq}\frac{\partial^3 (g_0)^{ks}}{\partial
y_0^i\partial y_0^j\partial y_0^q}\frac{\partial f_l^2}{\partial
x_0^s}=0,~1\leq l \leq {\bf
b}.$}\\

Similar to Theorem 9 in [PTN], we have\\

\noindent{\bf Theorem 5.9} {\it $M_0\times_{f_1} M_1\times \cdots
\times_{f_{\bf b}} M_{\bf b}$ is a relatively isotropic mean
Landsberg manifold, then $(M,F)$ is a weakly Landsberg manifold.}\\

Similar to Theorem 10 in [PTN], we have\\

\noindent{\bf Theorem 5.10} {\it $M_0\times_{f_1} M_1\times \cdots
\times_{f_{\bf b}} M_{\bf b}$ is locally dually flat if and only if
$(M_0,F_0)$ is locally dually flat and $f_l=f_l(x_l)$ and
$(M_i,f_iF_i)$ is locally dually flat for $1\leq l \leq {\bf
b}.$}\\

By Lemma 2 and Theorem 4 in [HR] and (5.25), we get\\

\noindent{\bf Theorem 5.11} {\it $M_0\times_{f_1} M_1\times \cdots
\times_{f_{\bf b}} M_{\bf b}$ is a locally Minkowski manifold if and
only if $f_l=f_l(x_l)$ and $(M_0,F_0),~(M_i,f_iF_i)$ are locally
Minkowski manifolds for $1\leq l \leq {\bf
b}.$}\\

Similar to Theorem 1 in [PT], we have\\

\noindent{\bf Theorem 5.12} {\it $M_0\times_{f_1} M_1\times \cdots
\times_{f_{\bf b}} M_{\bf b}$ is $C$-reducible, then it is a
Riemannian manifold.}\\

Let $$R_{ab}^c=\frac{\delta G^c_a}{\delta x^b}-\frac{\delta
G^c_b}{\delta x^a},\eqno(5.48)$$
$$\delta y_0^i=dy_0^i+G^i_jdx_0^j+\sum_{l=1}^{{\bf
b}}G^i_{\overline{l-1}+j}dx^j_l\eqno(5.49)$$ $$\delta
y_r^i=dy_l^{i}+G^{\overline{r-1}+i}_jdx_0^j+\sum_{l=1}^{{\bf
b}}G^{\overline{r-1}+i}_{\overline{l-1}+j}dx^j_l.\eqno(5.50)$$ Then
the multiply twisted Miron metric on $TM^0$ can be introduced as
follows:
$${\bf G}=(g_0)_{ij}dx_0^i\otimes dx_0^j+\sum_{l=1}^{{\bf
b}}f_l^2(g_l)_{ij}dx_l^i\otimes dx_l^j+\phi(F^2)(g_0)_{ij}\delta
y_0^i\otimes \delta y_0^j+\phi(F^2)\sum_{l=1}^{{\bf
b}}f_l^2(g_l)_{ij}\delta y_l^i\otimes \delta y_l^j.\eqno(5.51)$$\\

\noindent {\bf Proposition 5.13}{\it  Let $M=M_0\times_{f_1}
M_1\times \cdots \times_{f_{\bf b}} M_{\bf b}$, then the Levi-Civita
connection $\nabla$ on the Riemannian manifold $(TM^0,{\bf G})$ is
locally expressed as follows:}
$$\nabla_{\frac{\delta}{\delta x_0^i}}\frac{\delta}{\delta x_0^j}
=F_{ij}^s\frac{\delta}{\delta
x_0^s}+\left(\frac{1}{2}R_{ij}^s-\frac{1}{\phi(F^2)}(C_0)_{ij}^s\right)\frac{\partial}{\partial
y_0^s}+F_{ij}^{\overline{r-1}+t}{\frac{\delta}{\delta
x_r^t}}+\frac{1}{2}R_{ij}^{\overline{r-1}+t}\frac{\partial}{\partial
y_r^t},\eqno(5.52)$$

$$\nabla_{\frac{\delta}{\delta x_0^i}}\frac{\partial}{\partial
y_0^j}=\left[(C_0)_{ij}^s+\frac{1}{2}(g_0)^{sk}(g_0)_{cj}R^c_{ki}\phi(F^2)\right]\frac{\delta}{\delta
x_0^s}+F_{ij}^s\frac{\partial}{\partial y_0^s}$$
$$+\frac{1}{2f_r^2}(g^r)^{st}\left[f_r^2G_{ij}^{\overline{r-1}+l}(g_r)_{ls}-G_{i,\overline{r-1}+s}^l(g_0)_{lj}\right]
\frac{\partial}{\partial
y_r^s}+\frac{1}{2f_r^2}(g_r)^{st}(g_0)_{kj}R^k_{\overline{r-1}+s,i}\phi(F^2){\frac{\delta}{\delta
x_r^s}},\eqno(5.53)$$

$$\nabla_{\frac{\partial}{\partial
y_0^j}}\frac{\delta}{\delta x_0^i}=\nabla_{\frac{\delta}{\delta
x_0^i}}\frac{\partial}{\partial
y_0^j}-G_{ij}^l\frac{\partial}{\partial
y_0^l}-G_{ij}^{\overline{r-1}+l}\frac{\partial}{\partial
y_r^l},\eqno(5.54)$$

$$\nabla_{\frac{\delta}{\delta x_0^i}}\frac{\delta}{\delta
x_r^j}=\frac{1}{2}R_{i,\overline{r-1}+j}^s\frac{\partial}{\partial
y_0^s}+F_{i,\overline{r-1}+j}^s\frac{\delta}{\delta x_0^s}+
\frac{1}{2}R_{i,\overline{r-1}+j}^{\overline{q-1}+t}\frac{\partial}{\partial
y_q^t}+F_{i,\overline{r-1}+j}^{\overline{q-1}+t}\frac{\delta}{\delta
x_r^t},\eqno(5.55)$$
$$\nabla_{\frac{\delta}{\delta
x_r^j}}\frac{\delta}{\delta x_0^i}=
-\frac{1}{2}R_{i,\overline{r-1}+j}^s\frac{\partial}{\partial
y_0^s}+F_{i,\overline{r-1}+j}^s\frac{\delta}{\delta x_0^s}
-\frac{1}{2}R_{i,\overline{r-1}+j}^{\overline{q-1}+t}\frac{\partial}{\partial
y_q^t}+F_{i,\overline{r-1}+j}^{\overline{q-1}+t}\frac{\delta}{\delta
x_r^t},\eqno(5.56)$$

$$\nabla_{\frac{\partial}{\partial
y_r^i}}\frac{\partial}{\partial
y_t^j}=\frac{1}{2\phi(F^2)}\frac{\partial \phi(F^2)}{\partial
y_r^i}\frac{ \partial}{\partial y_t^j} +
\frac{1}{2\phi(F^2)}\frac{\partial \phi(F^2)}{\partial y_t^j}\frac{
\partial}{\partial y_r^i}, ~{\rm for}~r\neq t,\eqno(5.57)$$

$$\nabla_{\frac{\partial}{\partial
y_r^i}}\frac{\partial}{\partial
y_r^j}=-f_r^2(g_r)_{ij}\frac{\phi'(F^2)}{\phi(F^2)}y_0^l\frac{\partial}{\partial
y_0^l}$$
$$+\frac{1}{2}(g_0)^{lk}\phi(F^2)\left[-\frac{\delta}{\delta
x_0^k}(f_r^2(g_r)_{ij})+G_{k,\overline{r-1}+j}^{\overline{r-1}+s}f_r^2(g_r)_{si}+
G_{k,\overline{r-1}+i}^{\overline{r-1}+s}f_r^2(g_r)_{sj}\right]\frac{\delta}{\delta
x_0^l}$$
$$+\left[(C_q)_{ij}^\alpha\delta_r^q+\delta_r^q\delta_j^\alpha\frac{\phi'(F^2)}{\phi(F^2)}f_r^2(y_r)_i+
\delta_r^q\delta_i^\alpha\frac{\phi'(F^2)}{\phi(F^2)}f_r^2(y_r)_j-f_r^2\frac{\phi'(F^2)}{\phi(F^2)}(g_r)_{ij}y_q^\alpha\right]\frac{\partial}{\partial
y_q^\alpha}$$
 $$ +\frac{1}{2f_q^2}(g_q)^{lk}\left[
-\frac{\delta}{\delta
x_q^k}(f_q^2(g_q)_{ij})+G_{\overline{q-1}+k,\overline{r-1}+j}^{\overline{r-1}+\alpha}f_r^2(g_r)_{i\alpha}+
G_{\overline{q-1}+k,\overline{r-1}+i}^{\overline{r-1}+\alpha}f_r^2(g_r)_{j\alpha}\right]\phi(F^2)\frac{\delta}{\delta
x_q^l},\eqno(5.58)$$

$$\nabla_{\frac{\delta}{\delta
x_0^i}}\frac{\partial}{\partial y_r^j}=\frac{1}{2}(g_0)^{\alpha
k}\left[G_{i,\overline{r-1}+j}^l(g_0)_{lk}
-G_{ik}^{\overline{r-1}+l}(g_r)_{lj}f_r^2\right]\frac{\partial}{\partial
y_0^\alpha}$$
$$+\frac{1}{2}(g_0)^{sk}R_{ki}^{\overline{r-1}+l}(g_r)_{lj}f_r^2\phi(F^2)\frac{\delta}{\delta
x_0^s}$$ $$+\frac{1}{2f_r^2}(g_r)^{\alpha s}
\left[\frac{\delta}{\delta
x_0^i}(f_r^2(g_r)_{js})+G_{i,\overline{r-1}+j}^{\overline{r-1}+l}f_r^2(g_r)_{ls}-
G_{i,\overline{r-1}+s}^{\overline{r-1}+l}f_r^2(g_r)_{lj}\right]
\frac{\partial}{\partial y_r^\alpha}$$
$$
+\frac{1}{2f_k^2}(g_k)^{ts}R_{\overline{r-1}+s,i}^l(g_0)_{lj}\phi(F^2)\frac{\delta}{\delta
x_k^t},\eqno(5.59)$$

$$\nabla_{\frac{\partial}{\partial y_r^j}}\frac{\delta}{\delta
x_0^i}=\nabla_{\frac{\delta}{\delta x_0^i}}\frac{\partial}{\partial
y_r^j}-G_{i,\overline{r-1}+j}^k\frac{\partial}{\partial
y_0^k}-G_{i,\overline{r-1}+j}^{\overline{t-1}+k}\frac{\partial}{\partial
y_t^k},\eqno(5.60)$$

$$\nabla_{\frac{\partial}{\partial y_0^i}}\frac{\partial}{\partial
y_0^j}=\left\{(C_0)_{ij}^l+\frac{\phi'(F^2)}{\phi(F^2)}\left[\delta_j^l(y_0)_i+\delta_i^l(y_0)_j-(g_0)_{ij}y^l\right]\right\}\frac{\partial}{\partial
y_0^l}$$
$$+\frac{1}{2}(g_0)^{\alpha k}\phi(F^2)\left[-\frac{\delta}{\delta
x_0^k}(g_0)_{ij}+G_{kj}^l(g_0)_{li}+G_{ki}^l(g_0)_{lj}\right]\frac{\delta}{\delta
x_0^\alpha}-(g_0)_{ij}\frac{\phi'(F^2)}{\phi(F^2)}y^t\frac{\partial}{\partial
y_r^t}$$
$$+\frac{1}{2f_r^2}(g_r)^{t k}\phi(F^2)\left[-\frac{\delta}{\delta
x_r^k}(g_0)_{ij}+G_{\overline{r-1}+k,j}^l(g_0)_{li}+G_{\overline{r-1}+k,i}^l(g_0)_{lj}\right]\frac{\delta}{\delta
x_r^t},\eqno(5.61)$$

$$\nabla_{\frac{\delta}{\delta
x_r^j}}\frac{\partial}{\partial y_0^i} =\frac{1}{2} (g_0)^{t
k}\left[\frac{\delta}{\delta x_r^j}(g_0)_{ik}
+G_{\overline{r-1}+j,i}^l(g_0)_{lk}-G_{\overline{r-1}+j,k}^l(g_0)_{li}\right]\frac{\partial}{\partial
y_0^t}$$
$$+\frac{1}{2}(g_0)^{tk}(g_0)_{li}\phi(F^2)R_{k,\overline{r-1}+j}^l\frac{\delta}{\delta
x_0^t} +\frac{1}{2f_s^2}(g_s)^{\beta k}
\left[G_{\overline{r-1}+j,i}^{\overline{s-1}+l}f_s^2(g_s)_{lk}-
G_{\overline{r-1}+j,\overline{s-1}+k}^{l}(g_0)_{li}\right]\frac{\partial}{\partial
y_s^\beta}$$ $$ + \frac{1}{2f_s^2}(g_s)^{\beta
k}(g_0)_{li}\phi(F^2)R_{\overline{s-1}+k,\overline{r-1}+j}^l\frac{\delta}{\delta
x_s^\beta},\eqno(5.62)$$

$$\nabla_{\frac{\partial}{\partial y_0^i}}\frac{\delta}{\delta
x_r^j}=\nabla_{\frac{\delta}{\delta x_r^j}}\frac{\partial}{\partial
y_0^i} -G_{\overline{r-1}+j,i}^l\frac{\partial}{\partial y_0^l}-
G_{\overline{r-1}+j,i}^{\overline{t-1}+l}\frac{\partial}{\partial
y_t^l},\eqno(5.63)$$

$$\nabla_{\frac{\partial}{\partial y_r^i}}\frac{\partial}{\partial
y_0^j}=\frac{\phi'(F^2)}{\phi(F^2)}(y_r)_if_r^2\frac{\partial}{\partial
y_0^j}+\frac{1}{2}(g_0)^{lk}\phi(F^2)\left[
G_{k,j}^{\overline{r-1}+l}f_r^2(g_r)_{li}+
G_{k,\overline{r-1}+i}^{l}(g_0)_{lj}\right]\frac{\delta}{\delta
x_0^l}$$
$$+\frac{\phi'(F^2)}{\phi(F^2)}(y_0)_j\frac{\partial}{\partial y_r^i}
+ \frac{1}{2f_t^2}(g_t)^{\beta k}\phi(F^2)\left[
G_{\overline{t-1}+k,j}^{\overline{r-1}+l}f_r^2(g_r)_{li}+
G_{\overline{t-1}+k,\overline{r-1}+i}^{l}(g_0)_{lj}\right]\frac{\delta}{\delta
x_t^\beta},\eqno(5.64)$$

$$\nabla_{\frac{\partial}{\partial
y_0^j}}\frac{\partial}{\partial
y_r^i}=\nabla_{\frac{\partial}{\partial
y_r^i}}\frac{\partial}{\partial y_0^j},\eqno(5.65)$$

$$\nabla_{\frac{\delta}{\delta
x_r^i}}\frac{\delta}{\delta
x_t^j}=\frac{1}{2}R_{\overline{r-1}+i,\overline{t-1}+j}^{\beta}\frac{\partial}{\partial
y_0^\beta}+F_{\overline{r-1}+i,\overline{t-1}+j}^{\beta}\frac{\delta}{\delta
x_0^\beta}$$
$$+\left(-\delta_r^t\delta_r^s\frac{1}{\phi(F^2)}(C_r)_{ij}^\beta+\frac{1}{2}
R_{\overline{r-1}+i,\overline{t-1}+j}^{\overline{s-1}+\beta}\right)\frac{\partial}{\partial
y_s^\beta}+F_{\overline{r-1}+i,\overline{t-1}+j}^{\overline{l-1}+\beta}\frac{\delta}{\delta
x_l^\beta},\eqno(5.66)$$

$$\nabla_{\frac{\delta}{\delta
x_r^i}}\frac{\partial}{\partial y_t^j}= \frac{1}{2}(g_0)^{\beta
k}\left[G_{\overline{r-1}+i,\overline{t-1}+j}^l(g_0)_{lk}
-G_{\overline{r-1}+i,k}^{\overline{t-1}+l}(g_t)_{lj}f_t^2\right]\frac{\partial}{\partial
y_0^\beta}$$

$$+\frac{1}{2}(g_0)^{\beta k}R_{k,\overline{r-1}+i}^{\overline{t-1}+l}(g_t)_{lj}f_t^2\phi(F^2)\frac{\delta}{\delta
x_0^\beta}$$

 $$+\frac{1}{2f_s^2}(g_s)^{\beta k}
\left[\frac{\delta}{\delta
x_r^i}(\delta_t^sf_t^2(g_t)_{jk})+G_{\overline{r-1}+i,\overline{t-1}+j}^{\overline{s-1}+l}f_s^2(g_s)_{lk}-
G_{\overline{r-1}+i,\overline{s-1}+k}^{\overline{t-1}+l}f_t^2(g_t)_{lj}\right]
\frac{\partial}{\partial y_s^\beta}$$
$$
+\left[C_{\overline{r-1}+i,\overline{t-1}+j}^{\overline{s-1}+\beta}+
\frac{1}{2f_s^2}f_t^2(g_s)^{\beta
k}R_{\overline{s-1}+k,\overline{r-1}+i}^{\overline{t-1}+l}(g_t)_{lj}\phi(F^2)\right]\frac{\delta}{\delta
x_s^\beta},\eqno(5.67)$$

$$\nabla_{\frac{\partial}{\partial y_t^j}}\frac{\delta}{\delta
x_r^i}=\nabla_{\frac{\delta}{\delta x_r^i}}\frac{\partial}{\partial
y_t^j}-G_{\overline{r-1}+i,\overline{t-1}+j}^{l}\frac{\partial}{\partial
y_0^l}
-G_{\overline{r-1}+i,\overline{t-1}+j}^{\overline{s-1}+l}\frac{\partial}{\partial
y_s^l}.\eqno(5.68)$$\\

Let the horizontal bundle $HTM^0$ and vertical bundle $VTM^0$ have
the basis $\frac{\delta}{\delta x_0^i},~\frac{\delta}{\delta
x_r^i}$, $\frac{\partial}{\partial y_0^j},~\frac{\partial}{\partial
y_t^j}.$ We say that vertical bundle $VTM^0$ is totally geodesic in
$TTM^0$ if $\nabla_{\frac{\partial}{\partial
y^a}}\frac{\partial}{\partial y^b}\in\Gamma(VTM^0).$ Similarly, the
horizontal bundle $HTM^0$ is totally geodesic in $TTM^0$ if
$\nabla_{\frac{\delta}{\delta x^a}}\frac{\delta}{\delta
x^b}\in\Gamma(HTM^0).$ Similar to Proposition 3, 4 in [PT], we
have\\

\noindent {\bf Proposition 5.14}{\it  Let $M=M_0\times_{f_1}
M_1\times \cdots \times_{f_{\bf b}} M_{\bf b}$, then $VTM^0$ is
totally geodesic if and only if $F_{ab}^c=G_{ab}^c.$}\\

\noindent {\bf Proposition 5.15}{\it  Let $M=M_0\times_{f_1}
M_1\times \cdots \times_{f_{\bf b}} M_{\bf b}$, then $HTM^0$ is
totally geodesic if and only if $(M_i,F_i)$ is Riemannian for $0\leq
i\leq {\bf b}$ and $R_{ab}^c=0$}\\

 \noindent {\bf Acknowledgement.} This work
was supported by Fok Ying Tong Education
Foundation No. 121003.\\

\noindent{\large \bf References}\\

\noindent[ARS]L. Al\'{i}as, A. Romero, M. S\'{a}nchez, Spacelike
hypersurfaces of constant mean curvature and Clabi-Bernstein type
problems, Tohoku Math. J. 49(1997) 337-345.\\
\noindent[As1]G. Asanov, Finslerian extensions of Schwaraschild
metric, Fortschr. Phys. 40(1992) 667-693.\\
\noindent[As2]G. Asanov, Finslerian metric functions over the
product $R\times M$ and their potential applications, Rep. Math.
Phys. 41(1998) 117-132.\\
\noindent[BO]R. Bishop, B. O'Neill, Manifolds of negative curvature,
Trans. Am. Math. Soc. 145(1969) 1-49.\\
\noindent[BGV]M. Brozos-V\'{a}zquez, E. Garc\'{i}a-R\'{i}o, R.
V\'{a}zquez-Lorenzo, Some remarks on locally conformally flat static
space-times, J. Math. Phys. 46\\
\noindent[CK]J. Choi, M. Kim, The index form on the multiply warped
spacetime, Bull. Korean Math. Soc. 41(2004) No.4 691-697.\\
\noindent[DD]F. Dobarro, E. Dozo, Scalar curvature and warped
products of Riemannian manifolds, Trans. Am. Math. Soc. 303(1987)
161-168.\\
\noindent[DU1]F. Dobarro, B. \"{U}nal, Curvature of multiply warped
products, J. Geom. Phys. 55(2005) 75-106.\\
\noindent[DU2]F. Dobarro, B. \"{U}nal, Characterizing Killing vector
fields of standard static space-time, J. Geom. Phys. 62(2012)
1070-1087.\\
\noindent[EJK]P. Ehrlich, Y. Jung, S. Kim, Constant scalar
curvatures on warped product manifolds, Tsukuba J. Math. 20(1996)
No.1 239-265.\\
\noindent[EK]P. Ehrlich, S. Kim, The index form of a warped product,
preprint.\\
\noindent[FGKU]M. Fern\'{a}ndez-L\'{o}pez, E. Garc\'{i}a-R\'{i}o, D.
Kupeli, B. \"{U}nal, A curvature condition for a twisted product to
be a warped product, Manu. math. 106(2001), 213-217.\\
\noindent[FS]J. Flores, M. S\'{a}nchez, Geodesic connectedness of
multiwarped spacetimes, J. Diff. Eqs. 186(1)(2002) 1-30.\\
\noindent[Ha]H. Hayden, Subspace of a space with torsion, Proc.
Lond. Math. Soc. 34(1932) 27-50.\\
\noindent[HR]A. Hushmandi, M. Rezaii, On warped product Finsler
spaces of Landsberg type, J. Math. Phys. 52\\
\noindent[KPV]L. Kozma, I. Peter, c. Varga, Warped product of
Finsler manifolds, Ann. Univ. Sci. Budapest 44(2001)\\
\noindent[PT]E. Peyghan, A. Tayebi, On doubly warped product Finsler
manifolds, Nonlinear Analysis: Real world Applications 13(2012)
1703-1720.\\
\noindent[PTN]E. Peyghan, A. Tayebi, B. Najafi, Doubly warped
product Finsler manifolds with some non-Riemannian curvature
properties, arXiv:1110.6826.\\
\noindent[Sa]M. S\'{a}nchez, On the geometry of generalized
Robertson-Walker spacetimes: curvature and Killing fields, J. Geom.
Phys. 31(1999) 1-15.\\
 \noindent[SO]S. Sular, C. \"{O}zgur, Warped products with a
semi-symmetric metric connection, Taiwanese J. Math. 15(2011) no.4
1701-1719.\\
\noindent[U1]B. \"{U}nal, Multiply warped products, J. Geom. Phys.
34(2000) 287-301.\\
 \noindent[U2]B. \"{U}nal, Doubly warped products,
Diff. Geom. Appl. 15(2001) 253-263.\\
 \noindent[Ya]K. Yano, On semi-symmetric metric connection,
Rev. Roumaine Math. Pures Appl. 15(1970) 1579-1586.\\

 \indent{  School of Mathematics and Statistics,
Northeast Normal University, Changchun Jilin, 130024, China }\\
\indent E-mail: {\it wangy581@nenu.edu.cn}\\

\end{document}